\documentclass[11pt]{article}
\usepackage{amsmath}
\usepackage{amssymb}

\usepackage{theorem}
\numberwithin{equation}{section}

\newtheorem{theorem}{Theorem}[section]

\newtheorem{corollary}[theorem]{Corollary}
\newtheorem{conjecture}[theorem]{Conjecture}
\newtheorem{lemma}[theorem]{Lemma}
\newtheorem{definition}{Definition}[section]

\newcommand{\I}{P_{\leq CK}}

\begin{document}
\title{Global well-posedness and scattering for the defocusing, $L^{2}$-critical, nonlinear Schr{\"o}dinger equation when $d \geq 3$}
\date{\today}
\author{Benjamin Dodson}
\maketitle

\noindent \textbf{Abstract:} In this paper we prove that the defocusing, $d$-dimensional mass critical nonlinear Schr{\"o}dinger initial value problem is globally well-posed and scattering for $u_{0} \in L^{2}(\mathbf{R}^{d})$ and $d \geq 3$. To do this, we will prove a frequency localized interaction Morawetz estimate similar to the estimate made in \cite{CKSTT4}. Since we are considering an $L^{2}$ - critical initial value problem we will localize to low frequencies.

\section{Introduction} The $d$-dimensional, $L^{2}$ critical nonlinear Schr{\"o}dinger initial value problem is given by

\begin{equation}\label{0.1}
\aligned
i u_{t} + \Delta u &= F(u), \\
u(0,x) &= u_{0} \in L^{2}(\mathbf{R}^{d}),
\endaligned
\end{equation}

\noindent where $F(u) = \mu |u|^{4/d} u$, $\mu = \pm 1$, $u(t) : \mathbf{R}^{d} \rightarrow \mathbf{C}$. When $\mu = +1$ $(\ref{0.1})$ is said to be defocusing and when $\mu = -1$ $(\ref{0.1})$ is said to be focusing. The term $L^{2}$ - critical refers to scaling. If $u(t,x)$ solves $(\ref{0.1})$ on $[0, T]$ with initial data $u(0,x) = u_{0}(x)$, then

\begin{equation}\label{0.1.1}
\lambda^{d/2} u(\lambda^{2} t, \lambda x)
\end{equation}

\noindent solves $(\ref{0.1})$ on $[0, \frac{T}{\lambda^{2}}]$ with initial data $\lambda^{d/2} u_{0}(\lambda x)$. The scaling preserves the $L^{2}(\mathbf{R}^{d})$ norm. 

\begin{equation}\label{0.1.2}
\| \lambda^{d/2} u_{0}(\lambda x) \|_{L_{x}^{2}(\mathbf{R}^{d})} = \| u_{0}(x) \|_{L_{x}^{2}(\mathbf{R}^{d})}.
\end{equation}

\noindent It was observed in \cite{CaWe} that the solution to $(\ref{0.1})$ conserves the quantities mass,

\begin{equation}\label{0.2}
M(u(t)) = \int |u(t,x)|^{2} dx = M(u(0)),
\end{equation}

\noindent and energy

\begin{equation}\label{0.2.1}
E(u(t)) = \frac{1}{2} \int |\nabla u(t,x)|^{2} dx + \frac{\mu d}{2(d + 2)} \int |u(t,x)|^{\frac{2d + 4}{d}} dx = E(u(0)).
\end{equation}

\noindent \textbf{Remark:} When $\mu = +1$ this quantity is positive definite.\vspace{5mm}

\noindent A solution to $(\ref{0.1})$ obeys Duhamel's formula.

\begin{definition}\label{d0.1}
 $u : I \times \mathbf{R}^{d} \rightarrow \mathbf{C}$, $I \subset \mathbf{R}$ is a solution to $(\ref{0.1})$ if for any compact $J \subset I$, $u \in C_{t}^{0} L_{x}^{2}(J \times \mathbf{R}^{d}) \cap L_{t,x}^{\frac{2(d + 2)}{d}}(J \times \mathbf{R}^{d})$, and for all $t, t_{0} \in I$,

\begin{equation}\label{0.3.1}
 u(t) = e^{i(t - t_{0}) \Delta} u(t_{0}) - i \int_{t_{0}}^{t} e^{i(t - \tau) \Delta} F(u(\tau)) d\tau.
\end{equation}

\end{definition}

\noindent The space $L_{t,x}^{\frac{2(d + 2)}{d}}(J \times \mathbf{R}^{d})$ arises from the Strichartz estimates. This norm is also invariant under the scaling $(\ref{0.1.1})$.

\begin{definition}\label{d0.2}
 A solution to $(\ref{0.1})$ defined on $I \subset \mathbf{R}$ blows up forward in time if there exists $t_{0} \in I$ such that 

\begin{equation}\label{0.3.2}
 \int_{t_{0}}^{\sup(I)} \int |u(t,x)|^{\frac{2(d + 2)}{d}} dx dt = \infty.
\end{equation}

\noindent $u$ blows up backward in time if there exists $t_{0} \in I$ such that

\begin{equation}\label{0.3.3}
 \int_{\inf(I)}^{t_{0}} \int |u(t,x)|^{\frac{2(d + 2)}{d}} dx dt = \infty.
\end{equation}

\end{definition}

\begin{definition}\label{d0.0.2}
A solution $u(t,x)$ to $(\ref{0.1})$ is said to scatter forward in time if there exists $u_{+} \in L^{2}(\mathbf{R}^{d})$ such that

\begin{equation}\label{0.2.4}
\lim_{t \rightarrow \infty} \| e^{it \Delta} u_{+} - u(t,x) \|_{L^{2}(\mathbf{R}^{d})} = 0.
\end{equation}

\noindent A solution is said to scatter backward in time if there exists $u_{-} \in L^{2}(\mathbf{R}^{d})$ such that

\begin{equation}\label{0.2.4}
\lim_{t \rightarrow -\infty} \| e^{it \Delta} u_{-} - u(t,x) \|_{L^{2}(\mathbf{R}^{d})} = 0.
\end{equation}
\end{definition}

\begin{theorem}\label{t0.0.1}
For any $d \geq 1$, there exists $\epsilon(d) > 0$ such that if $\| u_{0} \|_{L^{2}(\mathbf{R}^{d})} < \epsilon(d)$, then $(\ref{0.1})$ is globally well-posed and scatters both forward and backward in time.
\end{theorem}

\noindent \emph{Proof:} See \cite{CaWe}, \cite{CaWe1}. $\Box$\vspace{5mm}

\noindent We will recall the proof of this theorem in $\S 2$. \cite{CaWe}, \cite{CaWe1} also proved $(\ref{0.1})$ is locally well-posed for $u_{0} \in L_{x}^{2}(\mathbf{R}^{d})$on some interval $[0, T]$, where $T(u_{0})$ depends on the profile of the initial data, not just its size in $L^{2}(\mathbf{R}^{d})$.

\begin{theorem}\label{t0.0.0.1}
 Given $u_{0} \in L^{2}(\mathbf{R}^{d})$ and $t_{0} \in \mathbf{R}$, there exists a maximal lifespan solution $u$ to $(\ref{0.1})$ defined on $I \subset \mathbf{R}$ with $u(t_{0}) = u_{0}$. Moreover,\vspace{5mm}

1. $I$ is an open neighborhood of $t_{0}$.

2. If $\sup(I)$ or $\inf(I)$ is finite, then $u$ blows up in the corresponding time direction.

3. The map that takes initial data to the corresponding solution is uniformly continuous on compact time intervals for bounded sets of initial data.

4. If $\sup(I) = \infty$ and $u$ does not blow up forward in time, then $u$ scatters forward to a free solution. If $\inf(I) = -\infty$ and $u$ does not blow up backward in time, then $u$ scatters backward to a free solution.
\end{theorem}

\noindent \emph{Proof:} See \cite{CaWe}, \cite{CaWe1}. $\Box$\vspace{5mm}

\noindent In the focusing case thee are known counterexamples to $(\ref{0.1})$ globally well-posed and scattering for all $u_{0} \in L^{2}(\mathbf{R}^{d})$. In the defocusing case there are no known counterexamples to global well-posedness and scattering for $u_{0} \in L^{2}(\mathbf{R}^{d})$ of arbitrary size. Therefore, it has been conjectured,

\begin{conjecture}\label{c0.0.2}
For $d \geq 1$, the defocusing, mass critical nonlinear Schr{\"o}dinger initial value problem $(\ref{0.1})$, $\mu = +1$ is globally well-posed for $u_{0} \in L^{2}(\mathbf{R}^{d})$ and all solutions scatter to a free solution as $t \rightarrow \pm \infty$.
\end{conjecture}

\noindent This conjecture has been affirmed in the radial case.

\begin{theorem}\label{t0.1}
When $d = 2$, $\mu = +1$, $(\ref{0.1})$ is globally well-posed and scattering for $u_{0} \in L^{2}(\mathbf{R}^{2})$ radial.
\end{theorem}

\noindent \emph{Proof:} See \cite{KTV}.

\begin{theorem}\label{t0.1.1}
When $d \geq 3$, $\mu = +1$, $(\ref{0.1})$ is globally well-posed and scattering for $u_{0} \in L^{2}(\mathbf{R}^{d})$ radial.
\end{theorem}

\noindent \emph{Proof:} See \cite{TVZ2}, \cite{KVZ}.\vspace{5mm}

\noindent In this paper we remove the radial condition for the case when $d \geq 3$ and prove

\begin{theorem}\label{t0.2}
$(\ref{0.1})$ is globally well-posed and scattering for $u_{0} \in L^{2}(\mathbf{R}^{d})$, $d \geq 3$.
\end{theorem}

\noindent \textbf{Remark:} \cite{KTV} and \cite{KVZ} also proved global well-posedness and scattering for the focusing, mass-critical initial value problem

\begin{equation}\label{0.2.4.1}
\aligned
iu_{t} + \Delta u &= -|u|^{4/d} u,\\
u(0,x) &= u_{0},
\endaligned
\end{equation}

\noindent with radial data and mass less than the mass of the ground state. Many of the tools used in this paper to prove global well-posedness and scattering when $\mu = +1$ can also be applied to the focusing problem with mass below the mass of the ground state. So whenever possible we will prove theorems for $\mu = \pm 1$.\vspace{5mm}

\noindent \textbf{Outline of the Proof.} We prove this theorem via the concentration compactness method, a modification of the induction on energy method. The induction on energy method was introduced in \cite{B2} to prove global well-posedness and scattering for the defocusing energy-critical initial value problem in $\mathbf{R}^{3}$ for radial data. \cite{B2} proved that it sufficed to treat solutions to the energy critical problem that were localized in both space and frequency. See \cite{CKSTT4}, \cite{RV}, \cite{V}, and \cite{TerryTao} for more work on the defocusing, energy critical initial value problem. \vspace{5mm}

\noindent This induction on energy method lead the development of the concentration compactness method. This method uses a concentration compactness technique to isolate a minimal mass/energy blowup solution. \cite{KTV} and \cite{KVZ} used concentration compactness to prove theorems $\ref{t0.1}$ and $\ref{t0.1.1}$. Since $(\ref{0.1})$ is globally well-posed for small $\| u_{0} \|_{L^{2}(\mathbf{R}^{d})}$, if $(\ref{0.1})$ is not globally well-posed for all $u_{0} \in L^{2}(\mathbf{R}^{d})$, then there must be a minimum $\| u_{0} \|_{L^{2}(\mathbf{R}^{d})} = m_{0}$ where global well-posedness fails. \cite{TVZ1} showed that for conjecture $\ref{c0.0.2}$ to fail, there must exist a minimal mass blowup solution with a number of additional properties. We show that such a solution cannot occur, proving theorem $\ref{t0.2}$. See \cite{KM}, \cite{Keraani}, \cite{Keraani1} for more information on this method.

\begin{definition}\label{d0.0.2.1}
 A set is precompact in $L^{2}(\mathbf{R}^{d})$ if it has compact closure in $L^{2}(\mathbf{R}^{d})$.
\end{definition}

\begin{definition}\label{d0.0.2.2}
 A solution $u(t,x)$ is said to be almost periodic if there exists a group of symmetries $G$ of the equation such that $\{ u(t) \}/G$ is a precompact set.
\end{definition}

\begin{theorem}\label{t0.4}
Suppose conjecture $\ref{c0.0.2}$ fails. Then there exists a maximal lifespan solution $u$ on $I \subset \mathbf{R}$, $u$ blows up both forward and backward in time, and $u$ is almost periodic modulo the group $G = (0, \infty) \times \mathbf{R}^{d} \times \mathbf{R}^{d}$ which consists of scaling symmetries, translational symmetries, and Galilean symmetries. That is, for any $t \in I$,

\begin{equation}\label{0.2.4.2}
u(t,x) = \frac{1}{N(t)^{d/2}} e^{ix \cdot \xi(t)} Q_{t}(\frac{x - x(t)}{N(t)}),
\end{equation}

\noindent where $Q_{t}(x) \in K \subset L^{2}(\mathbf{R}^{d})$, $K$ is a precompact subset of $L^{2}(\mathbf{R}^{d})$.\vspace{5mm}

\noindent Additionally, $[0, \infty) \subset I$, $N(t) \leq 1$ on $[0, \infty)$, $N(0) = 1$, and

\begin{equation}\label{0.2.4.3}
 \int_{0}^{\infty} \int |u(t,x)|^{\frac{2(d + 2)}{d}} dx dt = \infty.
\end{equation}

 \end{theorem}

\noindent \emph{Proof:} See \cite{TVZ1} and section four of \cite{TVZ2}. $\Box$\vspace{5mm}

\noindent \textbf{Remark:} This is also true of a minimal mass blowup solution to the focusing problem $(\ref{0.2.4.1})$.\vspace{5mm}

\noindent \textbf{Remark:} From the Arzela-Ascoli theorem, a set $K \subset L^{2}(\mathbf{R}^{d})$ is precompact if and only if there exists a compactness modulus function, $C(\eta) < \infty$ for all $\eta > 0$ such that

\begin{equation}\label{0.2.4.4}
\int_{|x| \geq C(\eta)} |f(x)|^{2} dx + \int_{|\xi| \geq C(\eta)} |\hat{f}(\xi)|^{2} d\xi < \eta.
\end{equation}

\noindent To verify conjecture $\ref{c0.0.2}$ in the case $d \geq 3$ it suffices to consider two scenarios separately,

\begin{equation}\label{0.2.4.5}
 \int_{0}^{\infty} N(t)^{3} dt = \infty,
\end{equation}

\noindent and

\begin{equation}\label{0.2.4.6}
 \int_{0}^{\infty} N(t)^{3} dt < \infty.
\end{equation}

\noindent The main new ingredient of this paper is to prove a long-time Strichartz estimate. The proof of this estimate relies on the bilinear Strichartz estimates and an induction on frequency argument.

\begin{theorem}\label{t0.5}
Suppose $J \subset [0, \infty)$ is compact, $d \geq 3$, $u$ is a minimal mass blowup solution to $(\ref{0.1})$ for $\mu = \pm 1$, and $\int_{J} N(t)^{3} dt = K$. Then there exists a function $\rho(N)$, $\rho(N) \leq 1$, $\lim_{N \rightarrow \infty} \rho(N) = 0$, such that for $N \leq K$,

\begin{equation}\label{0.2.4.7}
 \| P_{|\xi - \xi(t)| > N} u \|_{L_{t}^{2} L_{x}^{\frac{2d}{d - 2}}(J \times \mathbf{R}^{d})} \lesssim_{m_{0}, d} \rho(N) (\frac{K}{N})^{1/2}.
\end{equation}

\end{theorem}

\noindent To preclude the scenario $\int_{0}^{\infty} N(t)^{3} dt = \infty$ we will rely on a frequency localized interaction Morawetz estimate. (See \cite{CKSTT4} for such an estimate in the energy-critical case. \cite{CKSTT4} dealt with the energy-critical equation, $u(t) \in \dot{H}^{1}$, and thus truncated to high frequencies). The interaction Morawetz estimates scale like $\int_{J} N(t)^{3} dt$, and in fact are bounded below by some constant times $\int_{J} N(t)^{3} dt$. Since we are truncating to low frequencies, our method is very similar to the almost Morawetz estimates that are often used in conjunction with the I-method. (See \cite{B1}, \cite{CKSTT1}, \cite{CKSTT2}, \cite{CKSTT3}, \cite{CR}, \cite{CGT}, \cite{D}, \cite{D1}, \cite{DPST}, and \cite{DPST1} for more information on the I-method.) The estimates $(\ref{0.2.4.7})$ enable us to control the errors that arise from frequency truncation and prove\vspace{5mm}

\begin{theorem}\label{t0.8}
If $\int_{J} N(t)^{3} dt = K$, and $C$ is a large constant, independent of $K$, then

\begin{equation}\label{0.3}
\int_{J} \int_{\mathbf{R}^{d} \times \mathbf{R}^{d}} (-\Delta \Delta |x - y|) |\I u(t,x)|^{2} |\I u(t,y)|^{2} dx dy dt \lesssim_{m_{0}, d} o(K).
\end{equation}
\end{theorem}

\noindent This leads to a contradiction in the case when $\int_{0}^{\infty} N(t)^{3} dt = \infty$.\vspace{5mm}

\noindent To deal with the case when $\int_{0}^{\infty} N(t)^{3} dt < \infty$, we use a method similar to the method used in \cite{KTV}, \cite{TVZ}, and \cite{KVZ}. Such a minimal mass blowup solution must possess additional regularity in particular $u(t) \in L_{t}^{\infty} \dot{H}_{x}^{s}([0, \infty) \times \mathbf{R}^{d})$ for $0 < s < 1 + 4/d$. Since $\int N(t)^{3} dt < \infty$, $N(t) \searrow 0$ as $t \rightarrow \infty$, this contradicts conservation of energy. We rely on theorem $\ref{t0.5}$ to prove this additional regularity.\vspace{5mm}

\noindent \textbf{Outline of the Paper:} In $\S 2$, we describe some harmonic analysis and properties of the linear Schr{\"o}dinger equation that will be needed later in the paper. In particular we discuss Strichartz estimates. Global well-posedness and scattering for small mass will be an easy consequence of these estimates. We discuss the movement of $\xi(t)$ and $N(t)$ for a minimal mass blowup solution in this section. We also quote bilinear Strichartz estimates and the fractional chain rule.\vspace{5mm}

\noindent In $\S 3$ we prove theorem $\ref{t0.5}$. We use these estimates in $\S 4$ to obtain the frequency localized interaction Morawetz estimate and in $\S 5$ to obtain additional regularity.\vspace{5mm}

\noindent \textbf{Acknowledgements:} I am grateful to Monica Visan for her helpful comments on a preliminary draft of this paper.

\section{The linear Schr{\"o}dinger equation}
In this section we will introduce some of the tools that will be needed later in the paper. \vspace{5mm}

\noindent \textbf{Linear Strichartz Estimates:}

\begin{definition}\label{d5.0}
A pair $(p,q)$ will be called an admissible pair for $d \geq 3$ if $\frac{2}{p} = d(\frac{1}{2} - \frac{1}{q})$, and $p \geq 2$.
\end{definition}

\begin{theorem}\label{t5.1}
If $u(t,x)$ solves the initial value problem

\begin{equation}\label{5.1}
\aligned
i u_{t} + \Delta u &= F(t), \\
u(0,x) &= u_{0},
\endaligned
\end{equation}

\noindent on an interval $I$, then

\begin{equation}\label{5.2}
\| u \|_{L_{t}^{p} L_{x}^{q}(I \times \mathbf{R}^{d})} \lesssim_{p,q,\tilde{p},\tilde{q}, d} \| u_{0} \|_{L^{2}(\mathbf{R}^{d})} + \| F \|_{L_{t}^{\tilde{p}'} L_{x}^{\tilde{q}'}(I \times \mathbf{R}^{d})},
\end{equation}

\noindent for all admissible pairs $(p,q)$, $(\tilde{p}, \tilde{q})$. $\tilde{p}'$ denotes the Lebesgue dual of $\tilde{p}$.
\end{theorem}

\noindent \emph{Proof:} See \cite{Tao} for the case when $p > 2$, $\tilde{p} > 2$, and \cite{KT} for the proof when $p = 2$, $\tilde{p} = 2$, or both. We will rely very heavily on the double endpoint case, or when both $p = 2$ and $\tilde{p} = 2$.\vspace{5mm}

\noindent We will also make heavy use of the bilinear Strichartz estimates throughout the paper.

\begin{lemma}\label{l5.2}
Suppose $\hat{v}(t,\xi)$ is supported on $|\xi - \xi_{0}| \leq M$ and $\hat{u}(t,\xi)$ is supported on $|\xi - \xi_{0}| > N$, $M << N$, $\xi_{0} \in \mathbf{R}^{d}$. Then, for the interval $I = [a,b]$, $d \geq 1$,

\begin{equation}\label{5.5}
\| uv \|_{L_{t,x}^{2}(I \times \mathbf{R}^{d})} \lesssim \frac{M^{(d - 1)/2}}{N^{1/2}} \| u \|_{S_{\ast}^{0}(I \times \mathbf{R}^{d})} \| v \|_{S_{\ast}^{0}(I \times \mathbf{R}^{d})},
\end{equation}

\noindent where

\begin{equation}\label{5.6}
\| u \|_{S_{\ast}^{0}(I \times \mathbf{R}^{d})} \equiv \| u(a) \|_{L^{2}(\mathbf{R}^{d})} + \| (i \partial_{t} + \Delta) u \|_{L_{t,x}^{\frac{2(d + 2)}{d + 4}}(I \times \mathbf{R}^{d})}.
\end{equation}
\end{lemma}

\noindent \emph{Proof:} See \cite{V}.\vspace{5mm}

\noindent We will also need the Littlewood-Paley partition of unity. Let $\phi \in C_{0}^{\infty}(\mathbf{R}^{d})$, radial, $0 \leq \phi \leq 1$,

\begin{equation}\label{5.6.1}
\phi(x) = \left\{
            \begin{array}{ll}
              1, & \hbox{$|x| \leq 1$;} \\
              0, & \hbox{$|x| > 2$.}
            \end{array}
          \right.
\end{equation}

\noindent Define the frequency truncation

\begin{equation}\label{5.6.2}
\mathcal F(P_{\leq N} u) = \phi(\frac{\xi}{N}) \hat{u}(\xi).
\end{equation}

\noindent Let $P_{> N} u = u - P_{\leq N} u$ and $P_{N} u = P_{\leq 2N} u - P_{\leq N} u$. For convenience of notation let $u_{N} = P_{N} u$, $u_{\leq N} = P_{\leq N} u$, and $u_{> N} = P_{> N} u$.\vspace{5mm}

\noindent The Strichartz estimates motivate the definition of the Strichartz space.

\begin{definition}\label{d5.1}
Define the norm

\begin{equation}\label{5.3}
\| u \|_{S^{0}(I \times \mathbf{R}^{d})} \equiv \sup_{(p,q) \text{ admissible }} \| u \|_{L_{t}^{p} L_{x}^{q}(I \times \mathbf{R}^{d})}.
\end{equation}

\begin{equation}\label{5.3.1}
S^{0}(I \times \mathbf{R}^{d}) = \{ u \in C_{t}^{0}(I, L^{2}(\mathbf{R}^{d})) : \| u \|_{S^{0}(I \times \mathbf{R}^{d})} < \infty \}.
\end{equation}

\noindent We also define the space $N^{0}(I \times \mathbf{R}^{d})$ to be the space dual to $S^{0}(I \times \mathbf{R}^{d})$ with appropriate norm. Then in fact,

\begin{equation}\label{5.4}
\| u \|_{S^{0}(I \times \mathbf{R}^{d})} \lesssim \| u_{0} \|_{L^{2}(\mathbf{R}^{d})} + \| F \|_{N^{0}(I \times \mathbf{R}^{d})}.
\end{equation}
\end{definition}

\begin{theorem}\label{t5.1.1}
$(\ref{0.1})$ is globally well-posed when $\| u_{0} \|_{L^{2}(\mathbf{R}^{d})}$ is small.
\end{theorem}

\noindent \emph{Proof:} By $(\ref{5.4})$ and the definition of $S^{0}$, $N^{0}$,

\begin{equation}\label{5.4.1}
\| u \|_{L_{t,x}^{\frac{2(d + 2)}{d}}((-\infty, \infty) \times \mathbf{R}^{d})} \lesssim_{d} \| u_{0} \|_{L^{2}(\mathbf{R}^{d})} + \| u \|_{L_{t,x}^{\frac{2(d + 2)}{d}}((-\infty, \infty) \times \mathbf{R}^{d})}^{1 + 4/d}.
\end{equation}

\noindent By the continuity method, if $\| u_{0} \|_{L^{2}(\mathbf{R}^{d})}$ is sufficiently small, then we have global well-posedness. We can also obtain scattering with this argument. $\Box$\vspace{5mm}

\noindent Now let

\begin{equation}\label{5.4.2}
A(m) = \sup \{ \| u \|_{L_{t,x}^{\frac{2(d + 2)}{d}}((-\infty, \infty) \times \mathbf{R}^{d})} : \text{ u solves $(\ref{0.1})$, } \| u(0) \|_{L^{2}(\mathbf{R}^{d})} = m \}.
\end{equation}

\noindent If we can prove $A(m) < \infty$ for any $m$, then we have proved global well-posedness and scattering. Indeed, partition $(-\infty, \infty)$ into a finite number of subintervals with $\| u \|_{L_{t,x}^{\frac{2(d + 2)}{d}}(I_{j} \times \mathbf{R}^{d})} \leq \epsilon$ for each subinterval and iterate the argument in the proof of theorem $\ref{t5.1.1}$.\vspace{5mm}

\noindent Using a stability lemma from \cite{TVZ1} we can prove that $A(m)$ is a continuous function of $m$, which proves that $\{ m : A(m) = \infty \}$ is a closed set. This implies that if global well-posedness and scattering does not hold in the defocusing case for all $u_{0} \in L^{2}(\mathbf{R}^{d})$, then there must be a minimum $m_{0}$ with $A(m_{0}) = \infty$. Furthermore, \cite{TVZ1} proved that for conjecture $\ref{c0.0.2}$ to fail, there must exist a maximal interval $I \subset \mathbf{R}$ with $\| u \|_{L_{t,x}^{\frac{2(d + 2)}{d}}(I \times \mathbf{R}^{d})} = \infty$, and $u$ blows up both forward and backward in time. Moreover, this minimal mass blowup solution must be concentrated in both space and frequency. For any $\eta > 0$, there exists $C(\eta) < \infty$ with

\begin{equation}\label{6.2}
\int_{|x - x(t)| \geq \frac{C(\eta)}{N(t)}} |u(t,x)|^{2} dx < \eta,
\end{equation}

\noindent and

\begin{equation}\label{6.3}
\int_{|\xi - \xi(t)| \geq C(\eta) N(t)} |\hat{u}(t,\xi)|^{2} d\xi < \eta.
\end{equation}

\noindent By the Arzela-Ascoli theorem this proves $\{ u(t,x) \} / G$ is a precompact. It is quite clear that shifting the origin generates a $d$-dimensional symmetry group for solutions to $(\ref{0.1})$, and by $(\ref{0.1.1})$ changing $N(t)$ by a fixed constant also generates the multiplicative symmetry group $(0, \infty)$ for solutions to $(\ref{0.1})$. The Galilean transformation generates the $d$-dimensional phase shift symmetry group.

\begin{theorem}\label{t6.2}
Suppose $u(t,x)$ solves

\begin{equation}\label{6.4}
\aligned
i u_{t} + \Delta u &= F(u), \\
u(0,x) &= u_{0}.
\endaligned
\end{equation}

\noindent Then $v(t,x) = e^{-it |\xi_{0}|^{2}} e^{ix \cdot \xi_{0}} u(t, x - 2 \xi_{0} t)$ solves the initial value problem

\begin{equation}\label{6.5}
\aligned
i v_{t} + \Delta v &= F(v), \\
v(0,x) &= e^{ix \cdot \xi_{0}} u(0,x).
\endaligned
\end{equation}

\end{theorem}

\noindent \emph{Proof:} This follows by direct calculation. $\Box$\vspace{5mm}

\noindent If $u(t,x)$ obeys $(\ref{6.2})$ and $(\ref{6.3})$ and $v(t,x) = e^{-it |\xi_{0}|^{2}} e^{ix \cdot \xi_{0}} u(t, x - 2 \xi_{0} t)$, then

\begin{equation}\label{6.7}
\int_{|\xi - \xi_{0} - \xi(t)| \geq C(\eta) N(t)} |\hat{v}(t,\xi)|^{2} d\xi < \eta,
\end{equation}

\begin{equation}\label{6.8}
\int_{|x - 2 \xi_{0} t - x(t)| \geq \frac{C(\eta)}{N(t)}} |v(t,x)|^{2} dx < \eta.
\end{equation}

\noindent \textbf{Remark:} This will be useful to us later because it shifts $\xi(t)$ by a fixed amount $\xi_{0} \in \mathbf{R}^{d}$. For example, this allows us to set $\xi(0) = 0$. We now need to obtain some information on the movement of $N(t)$ and $\xi(t)$.

\begin{lemma}\label{l1.1}
If $J$ is an interval with

\begin{equation}\label{1.1}
\| u \|_{L_{t,x}^{\frac{2(d + 2)}{d}}(J \times \mathbf{R}^{d})} \leq C,
\end{equation}

\noindent then for $t_{1}, t_{2} \in J$,

\begin{equation}\label{1.2}
N(t_{1}) \sim_{C, m_{0}} N(t_{2}).
\end{equation}
\end{lemma}

\noindent \emph{Proof:} See \cite{KTV}, corollary 3.6. $\Box$

\begin{lemma}\label{l1.2.1}
If $u(t,x)$ is a minimal mass blowup solution on an interval J,

\begin{equation}\label{1.12.1}
\int_{J} N(t)^{2} dt \lesssim \| u \|_{L_{t,x}^{\frac{2(d + 2)}{d}}(J \times \mathbf{R}^{d})}^{\frac{2(d + 2)}{d}} \lesssim 1 + \int_{J} N(t)^{2} dt.
\end{equation}
\end{lemma}

\noindent \emph{Proof:} See \cite{KVZ}.\vspace{5mm}

\begin{lemma}\label{l1.4}
 Suppose $u$ is a minimal mass blowup solution with $N(t) \leq 1$. Suppose also that $J$ is some interval partitioned into subintervals $J_{k}$ with $\| u \|_{L_{t,x}^{\frac{2(d + 2)}{d}}(J_{k} \times \mathbf{R}^{d})} = \epsilon$ on each $J_{k}$. Again let

\begin{equation}\label{1.17.1}
N(J_{k}) = \sup_{J_{k}} N(t).
\end{equation}

\noindent Then,

\begin{equation}\label{1.18}
\sum_{J_{k}} N(J_{k}) \sim \int_{J} N(t)^{3} dt.
\end{equation}
\end{lemma}

\noindent \emph{Proof:} Since $N(t_{1}) \sim N(t_{2})$ for $t_{1}, t_{2} \in J_{k}$ it suffices to show $|J_{k}| \sim \frac{1}{N(J_{k})^{2}}$. By Holder's inequality and $(\ref{6.2})$, $$(\frac{m_{0}}{2})^{\frac{2(d + 2)}{d}} \leq (\int_{|x - x(t)| \leq \frac{C(\frac{m_{0}^{2}}{1000})}{N(t)}} |u(t,x)|^{2} dx)^{\frac{d + 2}{d}} \lesssim_{m_{0}} \frac{1}{N(t)^{2}} \| u(t,x) \|_{L_{x}^{\frac{2(d + 2)}{d}}(\mathbf{R}^{d})}^{\frac{2(d + 2)}{d}}.$$

\noindent Therefore, $$\int_{J_{k}} N(t)^{2} dt \lesssim_{m_{0}} \epsilon,$$ so $|J_{k}| \lesssim \frac{1}{N(J_{k})^{2}}$. Moreover, by Duhamel's formula, if $\| u \|_{L_{t,x}^{\frac{2(d + 2)}{d}}(J_{k} \times \mathbf{R}^{d})} = \epsilon$ then $$\| e^{i(t - a_{k}) \Delta} u(a_{k}) \|_{L_{t,x}^{\frac{2(d + 2)}{d}}(J_{k} \times \mathbf{R}^{d})} \geq \frac{\epsilon}{2},$$ where $J_{k} = [a_{k}, b_{k}]$. By Sobolev embedding,

\begin{equation}\label{1.19}
 \| e^{i(t - a_{k}) \Delta} P_{|\xi - \xi(a_{k})| \leq C(\epsilon^{2}) N(a_{k})} u(a_{k}) \|_{L_{t,x}^{\frac{2(d + 2)}{d}}(J_{k} \times \mathbf{R}^{d})} \lesssim_{m_{0}} N(J_{k})^{2} |J_{k}|.
\end{equation}

\noindent Therefore, $|J_{k}| \gtrsim \frac{1}{N(J_{k})^{2}}$. Summing up over subintervals proves the lemma. $\Box$\vspace{5mm}

\noindent We can use this fact to control the movement of $\xi(t)$. This control is essential for the arguments in the paper.

\begin{lemma}\label{l1.2}
Partition $J = [0, T_{0}]$ into subintervals $J = \cup J_{k}$ such that

\begin{equation}\label{1.7}
\| u \|_{L_{t,x}^{\frac{2(d + 2)}{d}}(J_{k} \times \mathbf{R}^{d})} \leq \epsilon,
\end{equation}

\noindent where $\epsilon$ is the same $\epsilon$ as in lemma $\ref{l1.1}$. Let $N(J_{k}) = \sup_{t \in J_{k}} N(t)$. Then

\begin{equation}\label{1.8}
|\xi(0) - \xi(T_{0})| \lesssim \sum_{k} N(J_{k}),
\end{equation}

\noindent which is the sum over the intervals $J_{k}$.
\end{lemma}

\noindent \emph{Proof:} See lemma 5.18 of $\cite{KilVis}$. $\Box$\vspace{5mm}

\noindent Possibly after adjusting the modulus function $C(\eta)$ in $(\ref{6.2})$, $(\ref{6.3})$ by a constant, we can choose $\xi(t) : I \rightarrow \mathbf{R}^{d}$ such that

\begin{equation}\label{1.13}
 |\frac{d}{dt} \xi(t)| \lesssim_{d} N(t)^{3}.
\end{equation}

\noindent \textbf{Fractional Chain Rule:} Another essential tool that we will need is a good analysis of embedding Holder continuous functions into Sobolev spaces. Since $d \geq 3$ our analysis of $(\ref{0.1})$ will be complicated by the fact that the nonlinearity $F(u) = \mu |u|^{4/d} u$ is no longer algebraic. Because of this fact, the Fourier transform of $F(u)$ is not the convolution of Fourier transforms of $u$, and thus $F(P_{< N})$ need not be truncated in frequency. Instead, we will use the fractional chain rule.

\begin{lemma}\label{l5.4}
Let $G$ be a Holder continuous function of order $0 < \alpha < 1$. Then for every $0 < s < \alpha$, $1 < p < \infty$, $\frac{s}{\alpha} < \sigma < 1$,

\begin{equation}\label{5.9}
\| |\nabla|^{s} G(u) \|_{L_{x}^{p}(\mathbf{R}^{d})} \lesssim \| |u|^{\alpha - \frac{s}{\sigma}} \|_{L_{x}^{p_{1}}(\mathbf{R}^{d})} \| |\nabla|^{\sigma} u \|_{L_{x}^{\frac{s}{\sigma} p_{2}}(\mathbf{R}^{d})}^{s/\sigma}.
\end{equation}
\end{lemma}

\noindent \emph{Proof:} See \cite{V}.\vspace{5mm}

\begin{corollary}\label{c5.4.1}
Let $0 \leq s < 1 + 4/d$. Then on any spacetime slab $I \times \mathbf{R}^{d}$,

\begin{equation}\label{5.9.1}
\| |\nabla|^{s} F(u) \|_{L_{t,x}^{\frac{2(d + 2)}{d + 4}}(I \times \mathbf{R}^{d})} \lesssim \| |\nabla|^{s} u \|_{L_{t,x}^{\frac{2(d + 2)}{d}}(I \times \mathbf{R}^{d})} \| u \|_{L_{t,x}^{\frac{2(d + 2)}{d}}(I \times \mathbf{R}^{d})}^{4/d}.
\end{equation}
\end{corollary}

\noindent \emph{Proof:} See \cite{KVZ}.\vspace{5mm}

\begin{corollary}\label{c5.5}
For $0 \leq s < 1 + \frac{4}{d}$,

\begin{equation}\label{5.10}
\| |\nabla|^{s} F(u) \|_{L_{t}^{2} L_{x}^{\frac{2d}{d + 2}}(J \times \mathbf{R}^{d})} \lesssim \| u \|_{L_{t}^{\infty} L_{x}^{2}(J \times \mathbf{R}^{d})}^{4/d} \| |\nabla|^{s} u \|_{L_{t}^{2} L_{x}^{\frac{2d}{d - 2}}(J \times \mathbf{R}^{d})}.
\end{equation}
\end{corollary}

\noindent \emph{Proof:} We use an argument similar to the argument found in \cite{KVZ} to prove corollary $\ref{c5.4.1}$. The case $s \leq 1$ follows from $\nabla F(u) = O(|u|^{4/d}) (\nabla u)$ and interpolating with the estimate for $\mu |u|^{4/d} u$. Now consider $s > 1$.\vspace{5mm}

\noindent \textbf{Case 1:, $d = 4$}

\begin{equation}\label{5.10.1}
\aligned
&\| \Delta F(u) \|_{L_{t}^{2} L_{x}^{4/3}(J \times \mathbf{R}^{4})} \\= \| F_{z}(u) \Delta u + F_{\bar{z}}(u) \Delta \bar{u} + &F_{zz}(u) (\nabla u)^{2} + F_{\bar{z} \bar{z}} (\nabla \bar{u})^{2} + 2 F_{z\bar{z}}(u) |\nabla u|^{2} \|_{L_{t}^{2} L_{x}^{4/3}(J \times \mathbf{R}^{4})}.
\endaligned
\end{equation}

\noindent By interpolation $$\| \nabla u \|_{L_{t}^{4} L_{x}^{8/3}(J \times \mathbf{R}^{4})}^{2} \lesssim \| \Delta u \|_{L_{t}^{2} L_{x}^{4}(J \times \mathbf{R}^{4})} \| u \|_{L_{t}^{\infty} L_{x}^{2}(J \times \mathbf{R}^{4})},$$

\noindent which proves the corollary in this case.\vspace{5mm}

\noindent \textbf{Case 2: $d > 4$:} Use the chain rule and fractional product rule (see \cite{T2} for more details).

\begin{equation}\label{5.11}
\aligned
&\| |\nabla|^{s} F(u) \|_{L_{t}^{2} L_{x}^{\frac{2d}{d - 2}}(J \times \mathbf{R}^{d})} \lesssim \| F_{z}(u) + F_{\bar{z}}(u) \|_{L_{t}^{\infty} L_{x}^{d/2}(J \times \mathbf{R}^{d})} \| |\nabla|^{s} u \|_{L_{t}^{2} L_{x}^{\frac{2d}{d - 2}}(J \times \mathbf{R}^{d})} \\ &+ \| |\nabla|^{s - 1} [F_{z}(u) + F_{\bar{z}}(u)] \|_{L_{t}^{\frac{2s}{s - 1}} L_{x}^{q}(J \times \mathbf{R}^{d})} \| \nabla u \|_{L_{t}^{2s} L_{x}^{p}(J \times \mathbf{R}^{d})},
\endaligned
\end{equation}

\noindent with

\begin{align}
\frac{1}{p} &= \frac{(d - 2)}{2ds} + \frac{s - 1}{2s}, \\
\frac{1}{q} &= \frac{2}{d} + \frac{(s - 1)(d - 2)}{2ds} - \frac{s - 1}{2s}.
\end{align}

\noindent By interpolation,

\begin{equation}\label{5.12}
\| \nabla u \|_{L_{t}^{2s} L_{x}^{p}(J \times \mathbf{R}^{d})} \lesssim \| |\nabla|^{s} u \|_{L_{t}^{2} L_{x}^{\frac{2d}{d - 2}}(J \times \mathbf{R}^{d})}^{1/s} \| u \|_{L_{t}^{\infty} L_{x}^{2}(J \times \mathbf{R}^{d})}^{(s - 1)/s}.
\end{equation}

\noindent Now use lemma $\ref{l5.4}$. Choose $\sigma$ with $\frac{s - 1}{4/d} < \sigma < 1$. Let $\frac{1}{p_{1}} = \frac{2}{d} - \frac{s - 1}{2 \sigma}$ and $\frac{1}{p_{2}} = \frac{(s - 1)(d - 2)}{2ds} + \frac{(s - \sigma)(s - 1)}{2s \sigma}$. Both $F_{z}(z)$ and $F_{\bar{z}}(z)$ are Holder continuous functions of order $\frac{4}{d}$. Without loss of generality consider $F_{z}(u)$.

\begin{equation}\label{5.13}
\| |\nabla|^{s - 1} F_{z}(u(t)) \|_{L_{x}^{q}(\mathbf{R}^{d})} \lesssim \| |u(t)|^{4/d - \frac{s - 1}{\sigma}} \|_{L_{x}^{p_{1}}(\mathbf{R}^{d})} \| |\nabla|^{\sigma} u(t) \|_{L_{x}^{(\frac{s - 1}{\sigma}) p_{2}}(\mathbf{R}^{d})}^{\frac{s - 1}{\sigma}}.
\end{equation}

\noindent By interpolation

\begin{equation}\label{5.13.1}
\aligned
\| |\nabla|^{\sigma} u(t) \|_{L_{t}^{\frac{2s}{\sigma}} L_{x}^{(\frac{s - 1}{s}) p_{2}}(J \times \mathbf{R}^{d})}^{\frac{s - 1}{\sigma}} \lesssim \| |\nabla|^{s} u \|_{L_{t}^{2} L_{x}^{\frac{2d}{d - 2}}(J \times \mathbf{R}^{d})}^{\frac{s - 1}{s}} \| u \|_{L_{t}^{\infty} L_{x}^{2}(J \times \mathbf{R}^{d})}^{(\frac{s - 1}{\sigma}) (\frac{s - \sigma}{s})}.
\endaligned
\end{equation}

\noindent Finally,

\begin{equation}\label{5.14}
\| |u|^{4/d - \frac{s - 1}{\sigma}} \|_{L_{t}^{\infty} L_{x}^{p_{1}}(J \times \mathbf{R}^{d})} \lesssim \| u \|_{L_{t}^{\infty} L_{x}^{2}(J \times \mathbf{R}^{d})}^{2/p_{1}}.
\end{equation}

\noindent Summing up our terms, the corollary is proved in this case also.\vspace{5mm}

\noindent \textbf{Case 3, $d = 3$:} Take $2 \leq s < 7/3$.

\begin{equation}\label{5.15}
\aligned
&\| |\nabla|^{s} F(u) \|_{L_{t}^{2} L_{x}^{6/5}(J \times \mathbf{R}^{3})} \\ = \| |\nabla|^{s - 2} [F_{z}(u) \Delta u + F_{\bar{z}}(u) \Delta \bar{u} + &2 F_{z \bar{z}}(u) |\nabla u|^{2} + F_{zz}(u) (\nabla u)^{2} + F_{\bar{z} \bar{z}}(u) (\nabla \bar{u})^{2} \|_{L_{t}^{2} L_{x}^{6/5}(J \times \mathbf{R}^{3})}.
\endaligned
\end{equation}

\noindent $F_{zz}$, $F_{z \bar{z}}$, $F_{\bar{z} \bar{z}}$ are Holder continuous of order $1/3$, while $F_{z}$ and $F_{\bar{z}}$ are in fact differentiable, so use lemma $\ref{l5.4}$ and interpolate as in the previous case. $\Box$\vspace{5mm}

\noindent Finally, at various points in the proof of theorem $\ref{t0.2}$ we will also rely on the Sobolev embedding lemma.

\begin{lemma}\label{l5.3}
If $\frac{1}{p} = \frac{1}{2} - \frac{\rho}{d}$ and $\rho < \frac{d}{2}$, then $$\dot{H}^{\rho}(\mathbf{R}^{d}) \subset L^{p}(\mathbf{R}^{d}),$$ and $$\| u \|_{L^{p}(\mathbf{R}^{d})} \lesssim_{p,d} \| u \|_{\dot{H}^{\rho}(\mathbf{R}^{d})}.$$
\end{lemma}

\noindent We will also rely on the Hardy-Littlewood-Sobolev lemma.

\begin{lemma}\label{l5.4}
Suppose $\frac{r}{d} = 1 - (\frac{1}{p} - \frac{1}{q})$, $1 < p < \infty$, $1 < q < \infty$, and $0 < r < d$. Then let

\begin{equation}\label{5.16}
G(x) = \int \frac{1}{|x - y|^{r}} F(y) dy.
\end{equation}

\begin{equation}\label{5.17}
\| G \|_{L^{q}(\mathbf{R}^{d})} \lesssim \| F \|_{L^{p}(\mathbf{R}^{d})}.
\end{equation}

\noindent We will use this result in $\S 4$ a great deal.
\end{lemma}



\section{Long-time Strichartz Estimates}
\noindent In order to defeat the minimal mass blowup solution we will obtain Strichartz estimates over long time intervals. These estimates will be used in $\S 4$ to preclude scenario $(\ref{0.2.4.5})$ from occurring and in $\S 5$ to preclude scenario $(\ref{0.2.4.6})$.

\begin{theorem}\label{t2.3}
Suppose $u$ is a minimal mass blowup solution to $(\ref{0.1})$, $\mu = \pm 1$, $J$ is a compact interval with $N(t) \leq 1$, and

\begin{equation}\label{2.5}
\int_{J} N(t)^{3} dt = K.
\end{equation}

\noindent Then for $N \leq K$, there exists a constant $C_{3}(m_{0}, d)$ such that

\begin{equation}\label{2.7}
\| P_{|\xi - \xi(t)| > N} u \|_{L_{t}^{2} L_{x}^{\frac{2d}{d - 2}}(J \times \mathbf{R}^{d})} \leq C_{3}(m_{0}, d) \frac{K^{1/2}}{N^{1/2}}.
\end{equation}

\end{theorem}

\noindent \emph{Proof:} We prove this theorem by induction on $N$. Start with the base case.

\begin{lemma}\label{l2.5}
Since $J$ is compact and $N(t) \leq 1$, $$\| u \|_{L_{t,x}^{\frac{2(d + 2)}{d}}(J \times \mathbf{R}^{d})}^{\frac{2(d + 2)}{d}} = C(J) < \infty.$$ Therefore, theorem $\ref{t2.3}$ is true for $N \leq \frac{K}{C(J)}$.
\end{lemma}

\noindent \noindent \emph{Proof:} Partition $J$ into $\frac{C^{2 + 4/d}}{\epsilon^{2 + 4/d}}$ subintervals $J_{k}$ with $\| u \|_{L_{t,x}^{\frac{2(d + 2)}{d}}(J_{k} \times \mathbf{R}^{d})} = \epsilon$. By Duhamel's formula and Strichartz estimates,

\begin{equation}\label{2.8.1}
\| u \|_{S^{0}(J_{k} \times \mathbf{R}^{d})} \lesssim_{d} \| u_{0} \|_{L^{2}(\mathbf{R}^{d})} + \| u \|_{L_{t,x}^{\frac{2(d + 2)}{d}}(J_{k} \times \mathbf{R}^{d})}^{1 + 4/d} \lesssim_{m_{0}, d} 1,
\end{equation}

\noindent which implies

\begin{equation}\label{2.8.2}
\| u \|_{L_{t}^{2} L_{x}^{\frac{2d}{d - 2}}(J \times \mathbf{R}^{d})} \leq C_{1}(m_{0}, d) C(J)^{1/2}.
\end{equation}

\noindent This implies theorem $\ref{t5.1}$ is true for the interval $J$ when $N \leq \frac{K}{C(J)}$.\vspace{5mm}

\noindent Next, we will make the inductive step. In the interest of first exposing the main idea, we will obtain an estimate conducive to induction when $\xi(t) \equiv 0$. After this, we will treat the case when $\xi(t)$ is time dependent, which necessarily introduces a few additional complications.\vspace{5mm}

\noindent \textbf{Remark:} The case $\xi(t) \equiv 0$ is already fairly interesting on its own. It includes the radial case, but also includes the case when $u(0,x)$ is symmetric across the $x_{1}, ..., x_{d}$ axes.\vspace{5mm}

\begin{lemma}\label{l2.8}
If $\xi(t) \equiv 0$, then there exists a function $\delta(C_{0})$, $\delta(C_{0}) \rightarrow 0$ as $C_{0} \rightarrow \infty$, such that when $d = 3$,

\begin{equation}\label{2.71}
\aligned
\| u_{> N}  \|_{L_{t}^{2} L_{x}^{6}(J \times \mathbf{R}^{3})} &\lesssim_{m_{0}, d, s} \| u_{> N} \|_{L_{t}^{\infty} L_{x}^{2}(J \times \mathbf{R}^{3})} + \sum_{M \leq \eta N} (\frac{M}{N})^{s} \| u_{> M} \|_{L_{t}^{2} L_{x}^{6}(J \times \mathbf{R}^{3})} \\ &+  \delta(C_{0}) \| u_{> \eta N} \|_{L_{t}^{2} L_{x}^{6}(J \times \mathbf{R}^{3})}
+ \frac{C_{0}^{3/2} K^{1/2}}{(\eta N)^{1/2}} (\sup_{J_{k}} \| u_{> \eta N} \|_{S_{\ast}^{0}(J_{k} \times \mathbf{R}^{3})}).
\endaligned
\end{equation}

\noindent When $d \geq 4$,

\begin{equation}\label{2.72}
\aligned
\| u_{> N} \|_{L_{t}^{2} L_{x}^{\frac{2d}{d - 2}}(J \times \mathbf{R}^{3})} &\lesssim_{m_{0}, d, s} \| u_{> N} \|_{L_{t}^{\infty} L_{x}^{2}(J \times \mathbf{R}^{d})} + \sum_{M \leq \eta N} (\frac{M}{N})^{s} \| u_{> M} \|_{L_{t}^{2} L_{x}^{\frac{2d}{d - 2}}(J \times \mathbf{R}^{d})} \\  + \delta(C_{0}) \| u_{> \eta N} \|_{L_{t}^{2} L_{x}^{\frac{2d}{d - 2}}(J \times \mathbf{R}^{d})}
&+ \frac{C_{0}^{4 - 6/d} K^{2/d}}{(\eta N)^{2/d}} (\sup_{J_{k}} \| u_{> \eta N} \|_{S_{\ast}^{0}(J_{k} \times \mathbf{R}^{d})})^{4/d} \| u_{> \eta N} \|_{L_{t}^{2} L_{x}^{\frac{2d}{d - 2}}(J \times \mathbf{R}^{d})}^{1 - 4/d}.
\endaligned
\end{equation}
\end{lemma}

\noindent \emph{Proof:} Define a cutoff $\chi(t) \in C_{0}^{\infty}(\mathbf{R}^{d})$ in physical space,

\begin{equation}\label{2.9}
\chi(t, x) = \left\{
            \begin{array}{ll}
              1, & \hbox{$|x - x(t)| \leq \frac{C_{0}}{N(t)}$;} \\
              0, & \hbox{$|x - x(t)| > \frac{2C_{0}}{N(t)}$.}
            \end{array}
          \right.
\end{equation}

\noindent $C_{0}$ will be specified later.

\begin{equation}\label{2.10}
\aligned
\| P_{> N} (|u(\tau)|^{4/d} u(\tau)) \|_{L_{t}^{2} L_{x}^{\frac{2d}{d + 2}}(J \times \mathbf{R}^{d})} &\lesssim_{d} \| P_{> N} (|u_{\leq \eta N}|^{4/d} u_{\leq \eta N}) \|_{L_{t}^{2} L_{x}^{\frac{2d}{d + 2}}(J \times \mathbf{R}^{d})} \\
+ \| (u_{> \eta N} ) |u_{> C_{0} N(t)}|^{4/d} \|_{L_{t}^{2} L_{x}^{\frac{2d}{d + 2}}(J \times \mathbf{R}^{d})}  &+ \| (u_{> \eta N} ) |(1 - \chi(t)) u_{\leq C_{0} N(t)}|^{4/d} \|_{L_{t}^{2} L_{x}^{\frac{2d}{d + 2}}(J \times \mathbf{R}^{d})} \\
&+ \| (u_{> \eta N} ) |\chi(t) u_{\leq C_{0} N(t)}|^{4/d} \|_{L_{t}^{2} L_{x}^{\frac{2d}{d + 2}}(J \times \mathbf{R}^{d})}.
\endaligned
\end{equation}

\noindent By Bernstein's inequality and $(\ref{5.10})$, for any $0 \leq s < 1 + 4/d$,

\begin{align}
\| P_{> N} (|u_{\leq \eta N}|^{4/d} u_{\leq \eta N}) \|_{L_{t}^{2} L_{x}^{\frac{2d}{d + 2}}(J \times \mathbf{R}^{d})} &\lesssim_{d} \frac{1}{N^{s}} \| |\nabla|^{s} u_{\leq \eta N} \|_{L_{t}^{2} L_{x}^{\frac{2d}{d - 2}}(J \times \mathbf{R}^{d})} \| u \|_{L_{t}^{\infty} L_{x}^{2}(J \times \mathbf{R}^{d})}^{4/d} \\ &\lesssim_{m_{0}, d} \sum_{M \leq \eta N} (\frac{M}{N})^{s} \| u_{> M} \|_{L_{t}^{2} L_{x}^{\frac{2d}{d - 2}}(J \times \mathbf{R}^{d})}.
\end{align}

\noindent For the next two terms we use $(\ref{6.2})$ and $(\ref{6.3})$. Since mass is concentrated in both frequency and space, we can deal with the mass outside these balls perturbatively.

$$\| (u_{> \eta N}) |u_{> C_{0} N(t)}|^{4/d} \|_{L_{t}^{2} L_{x}^{\frac{2d}{d + 2}}(J \times \mathbf{R}^{d})} + \|  (u_{> \eta N}) |(1 - \chi(t)) u_{\leq C_{0} N(t)}|^{4/d} \|_{L_{t}^{2} L_{x}^{\frac{2d}{d + 2}}(J \times \mathbf{R}^{d})}$$

$$\leq \| u_{> \eta N} \|_{L_{t}^{2} L_{x}^{\frac{2d}{d - 2}}(J \times \mathbf{R}^{d})} [\| (1 - \chi(t)) u \|_{L_{t}^{\infty} L_{x}^{2}(J \times \mathbf{R}^{d})}^{4/d} + \| u_{> C_{0} N(t)} \|_{L_{t}^{\infty} L_{x}^{2}(J \times \mathbf{R}^{d})}^{4/d}]$$

$$\leq \delta(C_{0}) \| u_{> \eta N} \|_{L_{t}^{2} L_{x}^{\frac{2d}{d - 2}}(J \times \mathbf{R}^{d})},$$

\noindent with $\delta(C_{0}) \rightarrow 0$ as $C_{0} \rightarrow \infty$ (see $(\ref{6.2})$, $(\ref{6.3})$). Finally, take

\begin{equation}\label{2.11}
\| (P_{> \eta N} u) |\chi(t) u_{\leq C_{0} N(t)}|^{4/d} \|_{L_{t}^{2} L_{x}^{\frac{2d}{d - 2}}(J \times \mathbf{R}^{d})}.
\end{equation}

\noindent We will use $(\ref{5.5})$ to estimate $(\ref{2.11})$ on each subinterval $J_{k}$ and then sum over all the subintervals.

\begin{equation}\label{2.11.1}
\| u \|_{S_{\ast}^{0}(J_{k} \times \mathbf{R}^{d})} = \| u_{0} \|_{L^{2}(\mathbf{R}^{d})} + \| |u|^{4/d} u \|_{L_{t,x}^{\frac{2(d + 2)}{d + 4}}(J_{k} \times \mathbf{R}^{d})} \lesssim_{m_{0}, d} 1.
\end{equation}

\noindent \textbf{When d = 3:} Recall that $N(J_{k}) = \sup_{t \in J_{k}} N(t)$. Applying the bilinear estimates, mass conservation $\| u \|_{L_{t}^{\infty} L_{x}^{2}(J \times \mathbf{R}^{d})} = m_{0}$, and Holder's inequality,

$$\|  (P_{> \eta N} u) |\chi(t) u_{\leq C_{0} N(t)}|^{4/3} \|_{L_{t}^{2} L_{x}^{6/5}(J_{k} \times \mathbf{R}^{3})}$$ $$\leq \| (P_{> \eta N} u)(u_{\leq C_{0} N(J_{k})}) \|_{L_{t,x}^{2}(J_{k} \times \mathbf{R}^{3})} \| \chi(t) \|_{L_{t}^{\infty} L_{x}^{6}(J_{k} \times \mathbf{R}^{d})} \| u \|_{L_{t}^{\infty} L_{x}^{2}(J \times \mathbf{R}^{d})}^{1/3}$$
$$\lesssim_{m_{0}, d} \frac{C_{0} N(J_{k})}{(\eta N)^{1/2}} (\frac{C_{0}}{N(J_{k})})^{1/2} \| u_{> \eta N} \|_{S_{\ast}^{0}(J_{k} \times \mathbf{R}^{d})} \| u \|_{S_{\ast}^{0}(J_{k} \times \mathbf{R}^{d})}.$$

\noindent Summing over the subintervals $J_{k}$ and using lemma $\ref{l1.4}$,

$$\|  (P_{> \eta N} u) |\chi u_{\leq C_{0} N(t)}|^{4/3} \|_{L_{t}^{2} L_{x}^{6/5}(J \times \mathbf{R}^{3})} \lesssim_{m_{0}, d} \frac{C_{0}^{3/2}}{\eta^{1/2}} \frac{K^{1/2}}{N^{1/2}} (\sup_{J_{k}} \| u_{> \eta N} \|_{S_{\ast}^{0}(J_{k} \times \mathbf{R}^{d})}).$$

\noindent \textbf{When $d \geq 4$:}\vspace{5mm}

\noindent To simplify notation let $\frac{1}{q} = \frac{2(d - 2)}{d^{2}}$ and $\frac{1}{p} = \frac{1}{q} + \frac{2}{d}$.

$$\| (P_{> \eta N} u) |\chi(t) u_{\leq C_{0} N(t)}|^{4/d} \|_{L_{t}^{2} L_{x}^{\frac{2d}{d + 2}}(J \times \mathbf{R}^{d})}$$

$$\leq \| |(P_{> \eta N} u)(u_{\leq C_{0} N(t)})|^{4/d} (\chi(t))^{4/d} \|_{L_{t}^{d/2} L_{x}^{p}(J \times \mathbf{R}^{d})} \| (P_{> \eta N} u)^{1 - 4/d} \|_{L_{t}^{2d/(d - 4)} L_{x}^{\frac{2d^{2}}{(d - 2)(d - 4)}}(J \times \mathbf{R}^{d})}.$$

\noindent Now,

$$\| [(P_{> \eta N} u) (u_{\leq C_{0} N(t)})]^{4/d} (\chi(t))^{4/d} \|_{L_{t}^{d/2} L_{x}^{p}(J_{k} \times \mathbf{R}^{d})}$$

$$\leq \| (P_{> \eta N} u)(u_{\leq C_{0} N(t)}) \|_{L_{t,x}^{2}(J_{k} \times \mathbf{R}^{d})}^{4/d} \| (\chi(t))^{4/d} \|_{L_{t}^{\infty} L_{x}^{q}(J_{k} \times \mathbf{R}^{d})}$$

$$\lesssim_{d} \frac{(C_{0} N(J_{k}))^{\frac{2(d - 1)}{d}}}{(\eta N)^{2/d}} \| u_{> \eta N} \|_{S_{\ast}^{0}(J_{k} \times \mathbf{R}^{d})}^{4/d} \| u \|_{S_{\ast}^{0}(J_{k} \times \mathbf{R}^{d})}^{4/d} (\frac{C_{0}}{N(J_{k})})^{\frac{2(d - 2)}{d}}$$

$$\lesssim_{m_{0}, d} C_{0}^{4 - 6/d} (\frac{N(J_{k})}{\eta N})^{2/d} \| u_{> \eta N} \|_{S_{\ast}^{0}(J_{k} \times \mathbf{R}^{d})}^{4/d}.$$

\noindent Again summing over all subintervals,

$$\| [(P_{> \eta N} u)(u_{\leq C_{0} N(t)})]^{4/d} (\chi(t))^{4/d} \|_{L_{t}^{d/2} L_{x}^{p}(J \times \mathbf{R}^{d})}$$

$$\lesssim_{m_{0}, d}  (\sum N(J_{k}))^{2/d} \frac{C_{0}^{4 - 6/d}}{(\eta N)^{2/d}} (\sup_{J_{k}} \| u_{> \eta N} \|_{S_{\ast}^{0}(J_{k} \times \mathbf{R}^{d})})^{4/d}$$ $$\lesssim_{m_{0}, d} \frac{K^{2/d}}{N^{2/d}} \frac{C_{0}^{4 - 6/d}}{\eta^{2/d}} (\sup_{J_{k}} \| u_{> \eta N} \|_{S_{\ast}^{0}(J_{k} \times \mathbf{R}^{d})})^{4/d}.$$

\noindent Therefore,

$$\| (u_{> \eta N} ) |\chi(t) u_{\leq C_{0} N(t)})|^{4/d}  \|_{L_{t}^{2} L_{x}^{\frac{2d}{d + 2}}(J \times \mathbf{R}^{d})}$$ $$\lesssim_{m_{0}, d} \frac{C_{0}^{2 - 4/d}}{\eta^{2/d}} \frac{K^{2/d}}{N^{2/d}} (\sup_{J_{k}} \| u_{> \eta N} \|_{S_{\ast}^{0}(J_{k} \times \mathbf{R}^{d})})^{4/d} \| u_{> \eta N} \|_{L_{t}^{2} L_{x}^{\frac{2d}{d - 2}}(J \times \mathbf{R}^{d})}^{1 - 4/d}.$$

\noindent By Strichartz estimates, when $d = 3$,

\begin{equation}\label{2.12}
\aligned
\| u_{> N}  \|_{L_{t}^{2} L_{x}^{6}(J \times \mathbf{R}^{3})} \lesssim_{m_{0}, d, s} \| u_{> N}  \|_{L_{t}^{\infty} L^{2}(J \times \mathbf{R}^{3})} + \sum_{M \leq \eta N} (\frac{M}{N})^{s} \| u_{> M} \|_{L_{t}^{2} L_{x}^{6}(J \times \mathbf{R}^{3})} \\ +  \delta(C_{0}) \| u_{> \eta N} \|_{L_{t}^{2} L_{x}^{6}(J \times \mathbf{R}^{3})}
+ \frac{C_{0}^{3/2} K^{1/2}}{(\eta N)^{1/2}} (\sup_{J_{k}} \| u_{> \eta N} \|_{S_{\ast}^{0}(J_{k} \times \mathbf{R}^{3})})
\endaligned
\end{equation}

\noindent This proves lemma $\ref{l2.8}$ when $d = 3$. When $d \geq 4$,

\begin{equation}\label{2.13}
\aligned
\| u_{> N} \|_{L_{t}^{2} L_{x}^{\frac{2d}{d - 2}}(J \times \mathbf{R}^{d})} \lesssim_{m_{0}, d, s} \| u_{> N}  \|_{L_{t}^{\infty} L_{x}^{2}(J \times \mathbf{R}^{d})}
+ \sum_{M \leq \eta N} (\frac{M}{N})^{s} \| u_{> M} \|_{L_{t}^{2} L_{x}^{\frac{2d}{d - 2}}(J \times \mathbf{R}^{d})} \\ + \delta(C_{0}) \| u_{> \eta N} \|_{L_{t}^{2} L_{x}^{\frac{2d}{d - 2}}(J \times \mathbf{R}^{d})}
+ \frac{C_{0}^{4 - 6/d} K^{2/d}}{(\eta N)^{2/d}} (\sup_{J_{k}} \| u_{> \eta N} \|_{S_{\ast}^{0}(J_{k} \times \mathbf{R}^{d})})^{4/d} \| u_{> \eta N} \|_{L_{t}^{2} L_{x}^{\frac{2d}{d - 2}}(J \times \mathbf{R}^{d})}^{1 - 4/d}.
\endaligned
\end{equation}

\noindent This proves lemma $\ref{l2.8}$. $\Box$\vspace{5mm}

\noindent Formulas $(\ref{2.12})$ and $(\ref{2.13})$ are quite good enough for us to prove theorem $\ref{t2.3}$ by induction, as will be shown in a moment. When $\xi(t)$ is time dependent we will settle for a slightly more complicated estimate.\vspace{5mm}

\noindent \textbf{$\xi(t)$ time dependent:} When $\xi(t)$ is time dependent we run into a bit of difficulty with the projection of the Duhamel term. Consider the case when $J = [0, T]$, $d = 3$, $N(t) \equiv 1$ and $\xi(t) = (t, 0, 0)$ to illustrate this idea. The low frequencies at time $t = 0$ will be the high frequencies at some later time. Indeed, at time $t > N$, $\xi = 0$ will belong to the set $$\{ |\xi - \xi(t)| > N \}.$$ Therefore, we cannot use the exact same argument as in the case when $\xi(t) \equiv 0$ because the projection $$\| P_{|\xi - \xi(t)| > N} (|u|^{4/3} (u)) \|_{L_{t}^{2} L_{x}^{6/5}([0, T] \times \mathbf{R}^{n})}$$ cannot be controlled by $$\| P_{|\xi - \xi(t)| > \eta N} u \|_{L_{t}^{2} L_{x}^{6}([0, T] \times \mathbf{R}^{n})}.$$\vspace{5mm}

\noindent Instead, we will partition $J$ into subintervals where $|\xi(t_{1}) - \xi(t_{2})| \lesssim N$ on each of the subintervals and use the Duhamel formula on each subinterval separately. By lemma $\ref{l1.2}$, $$|\xi(a) - \xi(b)| \lesssim_{d} \int_{a}^{b} N(t)^{3} dt.$$ So if $\int_{a}^{b} N(t)^{3} dt << N$, we can use the Duhamel formula and the triangle inequality to say

\begin{align}
\| P_{|\xi - \xi(t)| > N} u \|_{L_{t}^{2} L_{x}^{6}([a, b] \times \mathbf{R}^{3})} &\lesssim_{d} \| P_{|\xi - \xi(a)| > \frac{N}{2}} u(a) \|_{L_{x}^{2}(\mathbf{R}^{3})} + \| P_{|\xi - \xi(a)| > \frac{N}{2}} (|u|^{4/3} u) \|_{L_{t}^{2} L_{x}^{6/5}(J \times \mathbf{R}^{3})} \\ &\lesssim_{d} \| P_{|\xi - \xi(a)| > \frac{N}{2}} u(a) \|_{L_{x}^{2}(\mathbf{R}^{3})} + \| P_{|\xi - \xi(\tau)| > \frac{N}{4}} (|u|^{4/3} u)(\tau) \|_{L_{t}^{2} L_{x}^{6/5}(J \times \mathbf{R}^{3})}.
\end{align}

\noindent The tradeoff is that we are required to compute $\| P_{|\xi - \xi(t)| > N} u \|_{L_{t}^{2} L_{x}^{6}}$ over a bunch of subsets of $J$ separately and then add up their $L_{t}^{2} L_{x}^{6}$ norms.

\begin{lemma}\label{l2.9}
\noindent Suppose $\xi(t)$ is time dependent, and $u$ satisfies the same conditions as theorem $\ref{t2.3}$.

\begin{equation}\label{2.47}
\| u_{|\xi - \xi(t)| \geq N} \|_{L_{t}^{2} L_{x}^{\frac{2d}{d - 2}}(J \times \mathbf{R}^{d})} \lesssim_{m_{0}, d, s} (\frac{K}{N} + 1)^{1/2} \| u_{|\xi - \xi(t)| \geq \frac{N}{2}} \|_{L_{t}^{\infty} L_{x}^{2}(J \times \mathbf{R}^{d})} + (\sharp B_{j})^{1/2}
\end{equation}

\begin{equation}\label{2.48}
+ \sum_{M \leq \eta N} (\frac{M}{N})^{s} \| u_{|\xi - \xi(t)| \geq M} \|_{L_{t}^{2} L_{x}^{\frac{2d}{d - 2}}(J \times \mathbf{R}^{d})} + \delta(C_{0}) \| u_{|\xi - \xi(t)| \geq \eta N} \|_{L_{t}^{2} L_{x}^{\frac{2d}{d - 2}}(J \times \mathbf{R}^{d})}
\end{equation}

\begin{equation}\label{2.49}
+ \left\{
  \begin{array}{ll}
    C_{0}^{3/2} (\frac{K}{\eta N})^{1/2} (\sup_{J_{k}} \| u_{|\xi - \xi(t)| \geq \eta N} \|_{S_{\ast}^{0}(J_{k} \times \mathbf{R}^{d})}), & \hbox{if $d = 3$;} \\
    C_{0}^{4 - 6/d} (\frac{K}{\eta N})^{2/d} \| u_{|\xi - \xi(t)| \geq \eta N} \|_{L_{t}^{2} L_{x}^{\frac{2d}{d - 2}}(J \times \mathbf{R}^{d})}^{1 - 4/d} (\sup_{J_{k}} \| u \|_{S_{\ast}^{0}(J_{k} \times \mathbf{R}^{d})})^{4/d}, & \hbox{if $d \geq 4$.}
  \end{array}
\right.
\end{equation}

\begin{equation}\label{2.50}
+ (\frac{K}{\eta N})^{1/2} \left\{
                             \begin{array}{ll}
                               \| u_{|\xi - \xi(t)| \geq \eta N} \|_{L_{t}^{\infty} L_{x}^{2}(J \times \mathbf{R}^{3})}, & \hbox{if $d = 3$;} \\
                               \| u_{|\xi - \xi(t)| \geq \eta N} \|_{L_{t}^{\infty} L_{x}^{2}(J \times \mathbf{R}^{3})}^{4/d}, & \hbox{if $d \geq 4$.}
                             \end{array}
                           \right.
\end{equation}

\noindent $(\sharp B_{j})$ is the number of subintervals $J_{k}$ with $\| u \|_{L_{t,x}^{\frac{2(d + 2)}{d}}(J_{k} \times \mathbf{R}^{d})} = \epsilon$ and $N(J_{k}) > \frac{\eta_{1} N}{2}$. As in the case when $\xi(t) \equiv 0$, $\delta(C_{0}) \rightarrow 0$ as $C_{0} \rightarrow \infty$.
\end{lemma}

\noindent \emph{Proof:} By lemma $\ref{l1.2}$ we can choose $\eta_{1}(d)$ sufficiently small so that $|\xi(t_{1}) - \xi(t_{2})| \leq \frac{N(J_{k})}{100 \eta_{1}(d)}$ for $t_{1}, t_{2} \in J_{k}$. Since $J$ is compact and $N(t) \leq 1$, $J$ is the union of a finite number of subintervals $J_{k}$ with $\| u \|_{L_{t,x}^{\frac{2(d + 2)}{d}}(J_{k} \times \mathbf{R}^{d})} = \epsilon$. We will call these subintervals with $\| u \|_{L_{t,x}^{\frac{2(d + 2)}{d}}(J_{k} \times \mathbf{R}^{d})} = \epsilon$ the $\epsilon$ - subintervals.\vspace{5mm}

\noindent We will call the $\epsilon$ - subintervals with $N(J_{k}) > \frac{\eta_{1} N}{2}$ the bad subintervals. Then we will rewrite $J = \cup G_{j} \cup B_{j}$, where $B_{j}$ are the bad $\epsilon$ - subintervals and $G_{j}$ are the collections of good $\epsilon$ - subintervals in between the bad subintervals. Because $\sum N(J_{k}) \sim_{d} K$, $$(\sharp B_{j}) \lesssim_{d} \frac{2K}{N \eta_{1}}.$$ Next, cut each $G_{j}$ into some subcollections of $\epsilon$ - subintervals $G_{j} = \cup_{l} G_{j,l}$ with

\begin{equation}\label{2.13.2}
\sum N(J_{k}) \leq \eta_{1} N
\end{equation}

\noindent on each $G_{j,l}$, and such that one of three things is true about each $G_{j,l}$:\vspace{5mm}

1.
\begin{equation}\label{2.20}
\frac{\eta_{1} N}{2} \leq \sum_{J_{k} : J_{k} \cap G_{j,l} \neq \emptyset} N(J_{k}) \leq \eta_{1} N,
\end{equation}

2. $G_{j,l}$ is adjacent to $B_{j + 1}$,

\noindent or

3. $G_{j,l}$ is at the end of $J$.\vspace{5mm}

\noindent It is always possible to do this, because if $G_{j,l}$ is not adjacent to $B_{j + 1}$ or the end of $J$, and $$\sum_{J_{k} : J_{k} \cap G_{j,l} \neq \emptyset} N(J_{k}) < \frac{\eta_{1} N}{2},$$ we can add the $\epsilon$ - subinterval adjacent to $G_{j,l}$ to $G_{j,l}$ and still have $$\sum_{J_{k} : G_{j,l}} N(J_{k}) \leq \eta_{1} N.$$

\noindent Therefore,

\begin{equation}\label{2.20.0}
(\sharp G_{j,l}) \lesssim_{d} (\sharp B_{j}) + 1 + \frac{2K}{N \eta_{1}}.
\end{equation}

\noindent For the interval $B_{j}$ we will be content to simply say

\begin{equation}\label{2.20.1}
\|u \|_{L_{t}^{2} L_{x}^{\frac{2d}{d - 2}}(B_{j} \times \mathbf{R}^{d})} \lesssim 1 + \| u \|_{S^{0}(B_{j} \times \mathbf{R}^{d})}^{1 + 4/d} \lesssim_{m_{0}, d} 1.
\end{equation}

\noindent Now take $G_{j,l} = [a_{jl}, b_{jl}]$. By $(\ref{2.13.2})$, $|\xi(a_{jl}) - \xi(t)| \leq \frac{N}{100}$ when $t \in G_{j,l}$. This will give us something that is pretty close to $(\ref{2.12})$ and $(\ref{2.13})$ on each individual $G_{j,l}$.

\begin{lemma}\label{l2.6}
For $G_{j,l} = [a_{jl}, b_{jl}]$,

\begin{equation}\label{2.20.1.1}
\aligned
\| P_{|\xi - \xi(t)| > N} u \|_{L_{t}^{2} L_{x}^{\frac{2d}{d - 2}}(G_{j,l} \times \mathbf{R}^{d})} &\lesssim_{m_{0}, d, s} \| P_{|\xi - \xi(a_{jl})| > \frac{N}{2}} u(a_{jl}) \|_{L_{x}^{2}(\mathbf{R}^{d})} \\
&+ \delta(C_{0}) \| P_{|\xi - \xi(t)| > \eta N} u \|_{L_{t}^{2} L_{x}^{\frac{2d}{d - 2}}(G_{j,l} \times \mathbf{R}^{d})} \\
&+ \| (u_{|\xi - \xi(t)| > \eta N}) |\chi(t) u_{|\xi - \xi(t)| \leq C_{0} N(t)}|^{4/d} \|_{L_{t}^{2} L_{x}^{\frac{2d}{d - 2}}(G_{j,l} \times \mathbf{R}^{d})} \\
&+ \sum_{M \leq \eta N} (\frac{M}{N})^{s} \| u_{|\xi - \xi(t)| > M} \|_{L_{t}^{2} L_{x}^{\frac{2d}{d - 2}}(G_{j,l} \times \mathbf{R}^{d})}.
\endaligned
\end{equation}
\end{lemma}

\noindent \emph{Proof:} By Duhamel's formula the solution on $G_{j,l}$ has the form

\begin{equation}\label{2.20.2}
u(t,x) = e^{i(t - a_{jl}) \Delta} u(a_{jl}) - i \int_{a_{jl}}^{t} e^{i(t - \tau) \Delta} |u(\tau)|^{4/d} u(\tau) d\tau.
\end{equation}

\noindent Because $|\xi(a_{j,l}) - \xi(t)| \leq \frac{N}{100}$,

\begin{equation}\label{2.21}
\| P_{|\xi - \xi(t)| > N} u \|_{L_{t}^{2} L_{x}^{\frac{2d}{d - 2}}(G_{j,l} \times \mathbf{R}^{d})} \leq \| P_{|\xi - \xi(a_{jl})| > \frac{N}{2}} u \|_{L_{t}^{2} L_{x}^{\frac{2d}{d - 2}}(G_{j,l} \times \mathbf{R}^{d})}
\end{equation}

\begin{equation}\label{2.22}
\lesssim_{d} \| P_{|\xi - \xi(a_{jl})| > \frac{N}{2}} u(a_{jl}) \|_{L_{x}^{2}(\mathbf{R}^{d})}
+ \| P_{|\xi - \xi(a_{jl})| > \frac{N}{2}} (|u|^{4/d} u) \|_{L_{t}^{2} L_{x}^{\frac{2d}{d + 2}}(G_{j,l} \times \mathbf{R}^{d})}.
\end{equation}

\noindent Turning to the Duhamel term,

\begin{equation}\label{2.24}
\| P_{|\xi - \xi(a_{jl})| > \frac{N}{2}} (|u|^{4/d} u) \|_{L_{t}^{2} L_{x}^{\frac{2d}{d + 2}}(G_{j,l} \times \mathbf{R}^{d})}
\end{equation}

\begin{equation}\label{2.25}
\lesssim_{d} \| P_{|\xi - \xi(a_{jl})| > \frac{N}{2}} (|u_{|\xi - \xi(t)| \leq \eta N} |^{4/d} u_{|\xi - \xi(t)| \leq \eta N}) \|_{L_{t}^{2} L_{x}^{\frac{2d}{d + 2}}(G_{j,l} \times \mathbf{R}^{d})}
\end{equation}

\begin{equation}\label{2.27}
+ \| (u_{|\xi - \xi(t)| > \eta N}) |u_{|\xi - \xi(t)| > C_{0} N(t)}|^{4/d} \|_{L_{t}^{2} L_{x}^{\frac{2d}{d + 2}}(G_{j,l} \times \mathbf{R}^{d})}
\end{equation}

\begin{equation}\label{2.28}
+ \| (u_{|\xi - \xi(t)| > \eta N}) |(1 - \chi(t)) u|^{4/d} \|_{L_{t}^{2} L_{x}^{\frac{2d}{d + 2}}(G_{j,l} \times \mathbf{R}^{d})}
\end{equation}

\begin{equation}\label{2.26}
+ \| (u_{|\xi - \xi(t)| > \eta N}) |\chi(t) u_{|\xi - \xi(t)| \leq C_{0} N(t)}|^{4/d} \|_{L_{t}^{2} L_{x}^{\frac{2d}{d + 2}}(G_{j,l} \times \mathbf{R}^{d})}.
\end{equation}

\noindent By $(\ref{6.2})$ and $(\ref{6.3})$,

\begin{equation}\label{2.18}
\aligned
&\| (u_{|\xi - \xi(t)| \geq \eta N}) |(1 - \chi(t)) u_{|\xi - \xi(t)| \leq C_{0} N(t)}|^{4/d} \|_{L_{t}^{2} L_{x}^{\frac{2d}{d + 2}}(G_{j,l} \times \mathbf{R}^{d})} \\ &+ \| (u_{|\xi - \xi(t)| \geq \eta N}) | u_{|\xi - \xi(t)| > C_{0} N(t)}|^{4/d} \|_{L_{t}^{2} L_{x}^{\frac{2d}{d + 2}}(G_{j,l} \times \mathbf{R}^{d})} \\ &\leq \delta(C_{0}) \| u_{|\xi - \xi(t)| \geq \eta N} \|_{L_{t}^{2} L_{x}^{\frac{2d}{d - 2}}(G_{j,l} \times \mathbf{R}^{d})}.
\endaligned
\end{equation}

\noindent This takes care of $(\ref{2.27})$ and $(\ref{2.28})$. Next take $(\ref{2.25})$.

\begin{equation}\label{2.29}
\| P_{|\xi - \xi(a_{jl})| > \frac{N}{2}} (|u_{|\xi - \xi(t)| \leq \eta N}|^{4/d} u_{|\xi - \xi(t)| \leq \eta N}) \|_{L_{t}^{2} L_{x}^{\frac{2d}{d + 2}}(G_{j,l} \times \mathbf{R}^{d})}
\end{equation}

\begin{equation}\label{2.30}
= \| P_{|\xi - \xi(a_{jl}) + \xi(t)| > \frac{N}{2}}  (e^{-ix \cdot \xi(t)} |u_{|\xi - \xi(t)| \leq \eta N}|^{4/d} u_{|\xi - \xi(t)| \leq \eta N}) \|_{L_{t}^{2} L_{x}^{\frac{2d}{d + 2}}(G_{j,l} \times \mathbf{R}^{d})}
\end{equation}

\begin{equation}\label{2.31}
= \| P_{|\xi - \xi(a_{jl}) + \xi(t)| > \frac{N}{2}}  ( |e^{-ix \cdot \xi(t)} u_{|\xi - \xi(t)| \leq \eta N}|^{4/d}  (e^{-ix \cdot \xi(t)} u_{|\xi - \xi(t)| \leq \eta N})) \|_{L_{t}^{2} L_{x}^{\frac{2d}{d + 2}}(G_{j,l} \times \mathbf{R}^{d})}.
\end{equation}

\noindent Because $|\xi(a_{jl}) - \xi(t)| \leq \frac{N}{100}$ on $G_{j,l}$,

\begin{equation}\label{2.32}
(\ref{2.31}) \leq \| P_{|\xi| > \frac{N}{4}}  ( |e^{-ix \cdot \xi(t)} u_{|\xi - \xi(t)| \leq \eta N}|^{4/d} (e^{-ix \cdot \xi(t)} u_{|\xi - \xi(t)| \leq \eta N})) \|_{L_{t}^{2} L_{x}^{\frac{2d}{d + 2}}(G_{j,l} \times \mathbf{R}^{d})}.
\end{equation}

\noindent By Bernstein's inequality,

\begin{equation}\label{2.33}
(\ref{2.32}) \lesssim_{d} \frac{1}{N^{s}} \| |\nabla|^{s}  ( |e^{-ix \cdot \xi(t)} u_{|\xi - \xi(t)| \leq \eta N}|^{4/d} (e^{-ix \cdot \xi(t)} u_{|\xi - \xi(t)| \leq \eta N})) \|_{L_{t}^{2} L_{x}^{\frac{2d}{d + 2}}(G_{j,l} \times \mathbf{R}^{d})}.
\end{equation}

\noindent By corollary $\ref{c5.5}$, for $0 \leq s < 1 + 4/d$,

\begin{equation}\label{2.34}
(\ref{2.33}) \lesssim_{m_{0}, d, s} \frac{1}{N^{s}} \| |\nabla|^{s} (e^{-ix \cdot \xi(t)} u_{|\xi - \xi(t)| \leq \eta N}) \|_{L_{t}^{2} L_{x}^{\frac{2d}{d - 2}}(G_{j,l} \times \mathbf{R}^{d})},
\end{equation}

\begin{equation}\label{2.35}
\lesssim_{m_{0}, d, s} \sum_{M \leq \eta N} (\frac{M}{N})^{s} \| u_{|\xi - \xi(t)| > M} \|_{L_{t}^{2} L_{x}^{\frac{2d}{d - 2}}(G_{j,l} \times \mathbf{R}^{d})}.
\end{equation}

\noindent This finishes the proof of lemma $\ref{l2.6}$. $\Box$\vspace{5mm}

\noindent Returning to the proof of lemma $\ref{l2.9}$, summing the estimates $(\ref{2.20.1.1})$ over all the $G_{j,l}$ intervals, and using the crude estimate $(\ref{2.20.1})$ on each $B_{j}$,

\begin{equation}\label{2.36}
\| u_{|\xi - \xi(t)| > N} \|_{L_{t}^{2} L_{x}^{\frac{2d}{d - 2}}(J \times \mathbf{R}^{d})} \lesssim_{m_{0}, d, s} (\sharp G_{j,l})^{1/2} \| u_{|\xi - \xi(t)| > \frac{N}{2}}  \|_{L_{t}^{\infty} L_{x}^{2}(J \times \mathbf{R}^{d})} + (\sharp B_{j})^{1/2}
\end{equation}

\begin{equation}\label{2.37}
+ \delta(C_{0}) \| u_{|\xi - \xi(t)| > \eta N}  \|_{L_{t}^{2} L_{x}^{\frac{2d}{d - 2}}(J \times \mathbf{R}^{d})}
\end{equation}

\begin{equation}\label{2.38}
+ \| (u_{|\xi - \xi(t)| > \eta N}) |\chi(t) u_{|\xi - \xi(t)| \leq C_{0} N(t)}|^{4/d} \|_{L_{t}^{2} L_{x}^{\frac{2d}{d - 2}}(J \times \mathbf{R}^{d})}
\end{equation}

\begin{equation}\label{2.39}
+ \sum_{M \leq \eta N} (\frac{M}{N})^{s} \| u_{|\xi - \xi(t)| > M} \|_{L_{t}^{2} L_{x}^{\frac{2d}{d - 2}}(J \times \mathbf{R}^{d})}
\end{equation}

\noindent This is almost in an acceptable form for our purposes. All that we have left to do is make a bilinear estimate of $(\ref{2.38})$. Take one of the $\epsilon$ - subintervals $J_{k} = [a_{k}, b_{k}]$.\vspace{5mm}

\noindent Suppose $d = 3$ and $N(J_{k}) \leq \eta_{1} \eta N$. We have $|\xi(t) - \xi(a_{k})| \leq \frac{N(J_{k})}{\eta_{1}}$ for all $t \in J_{k}$. In particular, $$\{ \xi : |\xi - \xi(t)| \leq C_{0} N(J_{k}) \} \subset \{ \xi : |\xi - \xi(a_{k})| \leq (C_{0} + \frac{1}{\eta_{1}(d)}) N(J_{k}) \}$$ and $$\{ |\xi - \xi(t)| \geq \eta N \} \subset \{ \xi : |\xi - \xi(a_{k})| \geq \frac{\eta N}{2} \}.$$ Therefore,

\begin{equation}\label{2.41}
\aligned
\| (u_{|\xi - \xi(t)| \geq \eta N}) |\chi(t) u_{|\xi - \xi(t)| \leq C_{0} N(t)}|^{4/3} \|_{L_{t}^{2} L_{x}^{6/5}(J_{k} \times \mathbf{R}^{3})} \\
\leq \| (u_{|\xi - \xi(t)| \geq \eta N}) (u_{|\xi - \xi(t)| \leq C_{0} N(J_{k})}) \|_{L_{t,x}^{2}(J_{k} \times \mathbf{R}^{3})} \| \chi(t) \|_{L_{t}^{\infty} L_{x}^{6}(J_{k} \times \mathbf{R}^{3})} \| u \|_{L_{t}^{\infty} L_{x}^{2}(J_{k} \times \mathbf{R}^{3})}^{1/3}  \\
\lesssim_{m_{0}, d} C_{0}^{1/2} \| (u_{|\xi - \xi(a_{k})| \geq \frac{\eta N}{2}}) (u_{|\xi - \xi(a_{k})| \leq (C_{0} + \frac{1}{\eta_{1}(d)}) N(J_{k})} \|_{L_{t,x}^{2}(J_{k} \times \mathbf{R}^{3})} \\
\lesssim_{m_{0}, d} \frac{C_{0}^{3/2}}{\eta^{1/2}} \frac{N(J_{k})^{1/2}}{N^{1/2}} ( \| u_{|\xi - \xi(a_{k})| \geq \eta N} \|_{S_{\ast}^{0}(J_{k} \times \mathbf{R}^{3})}),
\endaligned
\end{equation}

\noindent \textbf{Remark:} We take it for granted that $C_{0}$ is large, in particular $>> \frac{1}{\eta_{1}}$.\vspace{5mm}

\noindent If $N(J_{k}) \geq \eta \eta_{1} N$ we simply say that since $\| u \|_{L_{t}^{2} L_{x}^{6}(J_{k} \times \mathbf{R}^{3})} \lesssim_{m_{0}, d} 1$ and $\| u \|_{L_{t}^{\infty} L_{x}^{2}(J_{k} \times \mathbf{R}^{d})} = m_{0}$,

\begin{equation}\label{2.42}
\| (u_{|\xi - \xi(t)| \geq \eta N}) |\chi(t) u_{|\xi - \xi(t)| \leq C_{0} N(t)}|^{4/3} \|_{L_{t}^{2} L_{x}^{6/5}(J_{k} \times \mathbf{R}^{3})} \lesssim_{m_{0}, d} \| u_{|\xi - \xi(t)| \geq \eta N} \|_{L_{t}^{\infty} L_{x}^{2}(J \times \mathbf{R}^{d})}.
\end{equation}

\noindent Because $\sum N(J_{k}) \sim_{d} K$ there are $\lesssim_{d} \frac{K}{\eta \eta_{1} N}$ intervals with $N(J_{k}) \geq \eta \eta_{1} N$.\vspace{5mm}

\noindent Now take $d \geq 4$. Let $\frac{1}{q} = \frac{2(d - 2)}{d^{2}}$ and $\frac{1}{p} = \frac{1}{q} + \frac{2}{d}$. If $N(J_{k}) \leq \eta_{1} \eta N$,

\begin{equation}\label{2.43}
\aligned
\| |(u_{|\xi - \xi(t)| \geq \eta N}) \chi(t) (u_{|\xi - \xi(t)| \leq C_{0} N(t)})|^{4/d} \|_{L_{t}^{d/2} L_{x}^{p}(J_{k} \times \mathbf{R}^{d})} \\ \leq \| (u_{|\xi - \xi(t)| \geq \eta N})(u_{|\xi - \xi(t)| \leq C_{0} N(t)}) \|_{L_{t,x}^{2}(J \times \mathbf{R}^{d})}^{4/d} \| \chi(t) \|_{L_{t}^{\infty} L_{x}^{q}(J \times \mathbf{R}^{d})} \\ 
\leq \| (u_{|\xi - \xi(a_{k})| \geq \frac{\eta N}{2}})(u_{|\xi - \xi(a_{k})| \leq (C_{0} + \frac{1}{\eta_{1}}) N(J_{k})}) \|_{L_{t,x}^{2}(J_{k} \times \mathbf{R}^{d})}^{4/d} \| \chi(t) \|_{L_{t}^{\infty} L_{x}^{q}(J \times \mathbf{R}^{d})} \\
\lesssim_{m_{0}, d} \frac{C_{0}^{4 - 6/d}}{\eta^{2/d}} \frac{N(J_{k})^{2/d}}{N^{2/d}}  (\| u_{|\xi - \xi(a_{k})| \geq \frac{\eta N}{2}} \|_{S_{\ast}^{0}(J_{k} \times \mathbf{R}^{d})})^{4/d}.
\endaligned
\end{equation}

\noindent If $N(J_{k}) \geq \eta \eta_{1} N$,

\begin{equation}\label{2.44}
\aligned
\| (u_{|\xi - \xi(t)| \geq \eta N}) |\chi(t) u_{|\xi - \xi(t)| \leq C_{0} N(t)}|^{4/d} \|_{L_{t}^{2} L_{x}^{\frac{2d}{d + 2}}(J_{k} \times \mathbf{R}^{d})} \lesssim_{m_{0}, d} \| u_{|\xi - \xi(t)| \geq \eta N} \|_{L_{t}^{\infty} L_{x}^{2}(J \times \mathbf{R}^{d})}^{4/d}.
\endaligned
\end{equation}

\noindent Once again there are $\lesssim_{d} \frac{K}{\eta_{1} \eta N}$ subintervals with $N(J_{k}) \geq \eta \eta_{1} N$.\vspace{5mm}

\noindent Therefore, if $d = 3$,

\begin{equation}\label{2.45}
\aligned
\| (u_{|\xi - \xi(t)| \geq \eta N}) |\chi(t) u_{|\xi - \xi(t)| \leq C_{0} N(t)}|^{4/3} \|_{L_{t}^{2} L_{x}^{6/5}(J \times \mathbf{R}^{3})} \\
 \lesssim_{m_{0}, d} \frac{K^{1/2} C_{0}^{3/2}}{(\eta N)^{1/2}} (\sup_{J_{k}; N(J_{k}) \leq \eta_{1} \eta N} \| u_{|\xi - \xi(a_{k})| \geq \frac{\eta N}{2}} \|_{S_{\ast}^{0}(J_{k} \times \mathbf{R}^{3})})  + \frac{K^{1/2}}{(\eta N)^{1/2}} \| u_{|\xi - \xi(t)| \geq \eta N} \|_{L_{t}^{\infty} L_{x}^{2}(J \times \mathbf{R}^{3})}.
\endaligned
\end{equation}

\noindent If $d \geq 4$,

\begin{equation}\label{2.46}
\aligned
&\| (u_{|\xi - \xi(t)| \geq \eta N}) |\chi(t) u_{|\xi - \xi(t)| \leq C_{0} N(t)}|^{4/d} \|_{L_{t}^{2} L_{x}^{\frac{2d}{d + 2}}(J \times \mathbf{R}^{d})} \\
 \lesssim_{m_{0}, d} \frac{K^{2/d} C_{0}^{4 - 6/d}}{(\eta N)^{2/d}} (\sup_{J_{k}; N(J_{k}) \leq \eta_{1} \eta N} &\| u_{|\xi - \xi(a_{k})| \geq \frac{\eta N}{2}} \|_{S_{\ast}^{0}(J_{k} \times \mathbf{R}^{d})}^{4/d}) \| u_{|\xi - \xi(t)| \geq \eta N} \|_{L_{t}^{2} L_{x}^{\frac{2d}{d - 2}}(J \times \mathbf{R}^{d})}^{1 - 4/d}  \\ &+ \frac{K^{1/2}}{(\eta N)^{1/2}} \| u_{|\xi - \xi(t)| \geq \eta N} \|_{L_{t}^{\infty} L_{x}^{2}(J \times \mathbf{R}^{d})}^{4/d}.
\endaligned
\end{equation}

\noindent Summing up $(\ref{2.36})$ - $(\ref{2.39})$ and substituting $(\ref{2.45})$ or $(\ref{2.46})$ for $(\ref{2.38})$, depending on dimension,

\begin{equation}\label{2.47}
\| u_{|\xi - \xi(t)| \geq N} \|_{L_{t}^{2} L_{x}^{\frac{2d}{d - 2}}(J \times \mathbf{R}^{d})} \lesssim_{m_{0}, d, s} (\frac{K}{N} + 1)^{1/2} \| u_{|\xi - \xi(t)| \geq \frac{N}{2}} \|_{L_{t}^{\infty} L_{x}^{2}(J \times \mathbf{R}^{d})} + (\sharp B_{j})^{1/2}
\end{equation}

\begin{equation}\label{2.48}
+ \sum_{M \leq \eta N} (\frac{M}{N})^{s} \| u_{|\xi - \xi(t)| \geq M} \|_{L_{t}^{2} L_{x}^{\frac{2d}{d - 2}}(J \times \mathbf{R}^{d})} + \delta(C_{0}) \| u_{|\xi - \xi(t)| \geq \eta N} \|_{L_{t}^{2} L_{x}^{\frac{2d}{d - 2}}(J \times \mathbf{R}^{d})}
\end{equation}

\begin{equation}\label{2.49}
+ \left\{
  \begin{array}{ll}
    C_{0}^{3/2} (\frac{K}{\eta N})^{1/2} (\sup_{J_{k}} \| u_{|\xi - \xi(a_{k})| \geq \frac{\eta N}{2}} \|_{S_{\ast}^{0}(J_{k} \times \mathbf{R}^{d})}), & \hbox{if $d = 3$;} \\
    C_{0}^{4 - 6/d} (\frac{K}{\eta N})^{2/d} \| u_{|\xi - \xi(a_{k})| \geq \frac{\eta N}{2}} \|_{L_{t}^{2} L_{x}^{\frac{2d}{d - 2}}(J \times \mathbf{R}^{d})}^{1 - 4/d} (\sup_{J_{k}} \| u \|_{S_{\ast}^{0}(J_{k} \times \mathbf{R}^{d})})^{4/d}, & \hbox{if $d \geq 4$.}
  \end{array}
\right.
\end{equation}

\begin{equation}\label{2.50}
+ (\frac{K}{\eta N})^{1/2} \left\{
                             \begin{array}{ll}
                               \| u_{|\xi - \xi(t)| \geq \eta N} \|_{L_{t}^{\infty} L_{x}^{2}(J \times \mathbf{R}^{3})}, & \hbox{if $d = 3$;} \\
                               \| u_{|\xi - \xi(t)| \geq \eta N} \|_{L_{t}^{\infty} L_{x}^{2}(J \times \mathbf{R}^{3})}^{4/d}, & \hbox{if $d \geq 4$.}
                             \end{array}
                           \right.
\end{equation}

\noindent We have used $(\sharp G_{j,l}) \lesssim_{d} \sharp (B_{j}) + 1 + (\frac{2K}{\eta_{1} N})$ and $\sharp B_{j} \lesssim_{d} \frac{2K}{\eta_{1} N}$ in $(\ref{2.47})$. The proof of lemma $\ref{l2.9}$ is now complete. $\Box$\vspace{5mm}

\noindent Now we are ready to prove theorem $\ref{t2.3}$. Let $s = 1$. For now make the crude estimates $\| u \|_{L_{t}^{\infty} L_{x}^{2}(J \times \mathbf{R}^{d})} \lesssim_{m_{0}} 1$ and $$\sup_{J_{k}} \| u_{|\xi - \xi(a_{k})| \geq \frac{\eta N}{2}} \|_{S_{\ast}^{0}(J_{k} \times \mathbf{R}^{d})} \leq \sup_{J_{k}} \| u \|_{S_{\ast}^{0}(J_{k} \times \mathbf{R}^{d})} \lesssim_{m_{0}, d} 1.$$ By $(\ref{2.47})$ - $(\ref{2.50})$,

\begin{align}
 \| u_{|\xi - \xi(t)| > N} \|_{L_{t}^{2} L_{x}^{6}(J \times \mathbf{R}^{3})} \leq C_{2}(m_{0}, d) (\frac{K}{N})^{1/2} + C_{2}(m_{0}, d) C_{0}^{3/2} (\frac{K}{\eta N})^{1/2} \\
+ C_{2}(m_{0}, d) \sum_{M \leq \eta N} (\frac{M}{N}) \| u_{|\xi - \xi(t)| > M} \|_{L_{t}^{2} L_{x}^{6}(J \times \mathbf{R}^{3})} + C_{2}(m_{0}, d) \delta(C_{0}) \| u_{|\xi - \xi(t)| > \eta N} \|_{L_{t}^{2} L_{x}^{6}(J \times \mathbf{R}^{3})}
\end{align}

\noindent We can prove theorem $\ref{t2.3}$ for $d = 3$ by induction. Suppose theorem $\ref{t2.3}$ is true for $M \leq \eta N$.

$$C_{2}(m_{0}, d) \sum_{M \leq \eta N} (\frac{M}{N}) \| u_{|\xi - \xi(t)| > M} \|_{L_{t}^{2} L_{x}^{6}(J \times \mathbf{R}^{3})} \leq 5 \eta^{1/2} C_{2}(m_{0}, d) C_{3}(m_{0}, d)  (\frac{K}{N})^{1/2}.$$

\noindent Choose $\eta(m_{0}, d)$ sufficiently small so that $\eta^{1/2} C_{2}(m_{0}, d) \leq \frac{1}{1000}$.\vspace{5mm}

\noindent Next,

$$\delta(C_{0}) C_{2}(m_{0}, d) \| u_{|\xi - \xi(t)| > \eta N} \|_{L_{t}^{2} L_{x}^{\frac{2d}{d - 2}}(J \times \mathbf{R}^{d})} \leq \delta(C_{0}) C_{2}(m_{0}, d) C_{3}(m_{0}, d) (\frac{K}{\eta N})^{1/2}.$$

\noindent Since $\delta(C_{0}) \rightarrow 0$ as $C_{0} \rightarrow \infty$, choose $C_{0}(\eta(m_{0}, d), m_{0}, d)$ sufficiently large so that $\delta(C_{0}) \frac{C_{2}(m_{0}, d)}{\eta^{1/2}} \leq \frac{1}{1000}$.\vspace{5mm}

\noindent Finally, choose $C_{3}(m_{0}, d)$ sufficiently large so that

$$C_{2}(m_{0}, d) + C_{2}(m_{0}, d) \frac{C_{0}(\eta(m_{0}, d), m_{0}, d)^{3/2}}{\eta(m_{0}, d)^{1/2}} \leq \frac{1}{1000} C_{3}(m_{0}, d).$$

\noindent This closes the induction and proves theorem $\ref{t2.3}$ when $d = 3$.\vspace{5mm}

\noindent We make a similar argument for $d \geq 4$.

\begin{align}
 \| u_{|\xi - \xi(t)| > N} \|_{L_{t}^{2} L_{x}^{\frac{2d}{d - 2}}(J \times \mathbf{R}^{d})} \leq C_{2}(m_{0}, d) (\frac{K}{N})^{1/2} + C_{2}(m_{0}, d) C_{0}^{4 - 6/d} (\frac{K}{\eta N})^{2/d} \| u_{|\xi - \xi(t)| > \eta N} \|_{L_{t}^{2} L_{x}^{\frac{2d}{d - 2}}(J \times \mathbf{R}^{d})}^{1 - 4/d}\\
+ C_{2}(m_{0}, d) \sum_{M \leq \eta N} (\frac{M}{N}) \| u_{|\xi - \xi(t)| > M} \|_{L_{t}^{2} L_{x}^{6}(J \times \mathbf{R}^{3})} + C_{2}(m_{0}, d) \delta(C_{0}) \| u_{|\xi - \xi(t)| > \eta N} \|_{L_{t}^{2} L_{x}^{\frac{2d}{d - 2}}(J \times \mathbf{R}^{3})}
\end{align}

\noindent Choose $\eta(m_{0}, d) > 0$ sufficiently small so that $\eta^{1/2} C_{2}(m_{0}, d) \leq \frac{1}{1000}$. Next, choose $C_{0}(\eta(m_{0}, d), m_{0}, d)$ sufficiently large so that $\delta(C_{0}) \frac{C_{2}(m_{0}, d)}{\eta^{1/2}} \leq \frac{1}{1000}$. Finally, choose $C_{3}(m_{0}, d)$ sufficiently large so that

$$C_{2}(m_{0}, d) + C_{2}(m_{0}, d) \frac{C_{0}^{4 - 6/d}}{\eta^{1/2}} \leq \frac{1}{1000} C_{3}(m_{0}, d)^{4/d}.$$ This closes the induction and proves theorem $\ref{t2.3}$ when $d \geq 4$. $\Box$\vspace{5mm}

\noindent For the upcoming section we will need $$\| u_{|\xi - \xi(t)| > N} \|_{L_{t}^{2} L_{x}^{\frac{2d}{d - 2}}(J \times \mathbf{R}^{d})}$$ to decay slightly faster than $(\frac{K}{N})^{1/2}$.

\begin{theorem}\label{t2.4}
There exists a function $\rho(N) \leq 1$,

\begin{equation}\label{2.51}
\lim_{N \rightarrow \infty} \rho(N) = 0,
\end{equation}

\noindent such that if $u$ is a minimal mass blowup solution to $(\ref{0.1})$, $\mu = \pm 1$ on the compact interval $J$ with $N(t) \leq 1$ and $\int_{J} N(t)^{3} dt = K$, then

\begin{equation}\label{2.52}
\| u_{|\xi - \xi(t)| > N} \|_{L_{t}^{2} L_{x}^{\frac{2d}{d - 2}}(J \times \mathbf{R}^{d})} \leq C_{3}(m_{0}, d) \rho(N) (\frac{K}{N})^{1/2}.
\end{equation}
\end{theorem}

\noindent \emph{Proof:} We will modify the argument of the proof of theorem $\ref{t2.3}$ slightly, taking advantage of the decay afforded by $(\ref{6.3})$,

\begin{equation}\label{2.52.1}
 \lim_{N \rightarrow \infty} \| u_{|\xi - \xi(t)| > N} \|_{L_{t}^{\infty} L_{x}^{2}(J \times \mathbf{R}^{d})} = 0.
\end{equation}

\begin{lemma}\label{l2.5}
Let $J_{k}$ be an interval with $\| u \|_{L_{t,x}^{\frac{2(d + 2)}{d}}(J_{k} \times \mathbf{R}^{d})} = \epsilon$, $N(J_{k}) \leq 1$, and let $u$, $J$ satisfy the hypotheses of theorem $\ref{t2.4}$. Then there exists a function $\sigma(N)$, $\sigma(N) \lesssim_{m_{0}, d} 1$, $\lim_{N \rightarrow \infty} \sigma(N) = 0$, such that

\begin{equation}\label{2.52.2}
 \sup_{J_{k} = [a_{k}, b_{k}] \subset J} \| u_{|\xi - \xi(a_{k})| > N} \|_{S_{\ast}^{0}(J_{k} \times \mathbf{R}^{d})} \leq \sigma(N).
\end{equation}

\end{lemma}

\noindent \emph{Proof:} Since $N(J_{k}) \leq 1$, $|\xi(t) - \xi(a_{k})| \leq \frac{1}{100 \eta_{1}(d)}$ on $J_{k}$. Take $N \geq \frac{1000}{\eta_{1}(d)}$. The lemma follows from Strichartz estimates for $N \leq \frac{1000}{\eta_{1}(d)}$.

\begin{equation}\label{2.55}
\aligned
\| u_{|\xi - \xi(a_{k})| > N}  \|_{S_{\ast}^{0}(J_{k} \times \mathbf{R}^{d})} &\leq \| P_{|\xi - \xi(t)| > \frac{N}{2}} u\|_{L_{t}^{\infty} L_{x}^{2}(J_{k} \times \mathbf{R}^{d})} \\ &+ \| P_{|\xi - \xi(t)| > \frac{N}{2}} (|u|^{4/d} u) \|_{L_{t,x}^{\frac{2(d + 2)}{d + 4}}(J_{k} \times \mathbf{R}^{d})}.
\endaligned
\end{equation}

\noindent By Bernstein's inequality, $\| u \|_{L_{t,x}^{\frac{2(d + 2)}{d}}(J_{k} \times \mathbf{R}^{d})} \leq \epsilon$,

$$\| P_{|\xi - \xi(t)| \geq \frac{N}{2}} (|u_{|\xi - \xi(t)| \leq \frac{N^{1/2}}{2}}|^{4/d} u_{|\xi - \xi(t)| \leq \frac{N^{1/2}}{2}}) \|_{L_{t,x}^{\frac{2(d + 2)}{d + 4}}(J_{k} \times \mathbf{R}^{d})}$$

\begin{equation}\label{2.56}
\aligned
\lesssim \frac{1}{N} \| \nabla e^{-ix \cdot \xi(t)}(|u_{|\xi - \xi(t)| \leq \frac{N^{1/2}}{2}}|^{4/d} u_{|\xi - \xi(t)| \leq \frac{N^{1/2}}{2}} \|_{L_{t,x}^{\frac{2(d + 2)}{d + 4}}(J_{k} \times \mathbf{R}^{d})} \\
\lesssim_{d} \frac{1}{N} \| \nabla (e^{-ix \cdot \xi(t)} u_{|\xi - \xi(t)| \leq \frac{N^{1/2}}{2}}) \|_{L_{t,x}^{\frac{2(d + 2)}{d}}(J_{k} \times \mathbf{R}^{d})} \| u \|_{L_{t,x}^{\frac{2(d + 2)}{d}}(J_{k} \times \mathbf{R}^{d})} \lesssim_{m_{0}, d} N^{-1/2}.
\endaligned
\end{equation}

\noindent Also,

\begin{equation}\label{2.57}
\| |u_{|\xi - \xi(t)| \geq \frac{N^{1/2}}{2}}| |u|^{4/d} \|_{L_{t,x}^{\frac{2(d + 2)}{d + 4}}(J_{k} \times \mathbf{R}^{d})} \leq \| u_{|\xi - \xi(t)| \geq \frac{N^{1/2}}{2}} \|_{L_{t}^{\infty} L_{x}^{2}(J_{k} \times \mathbf{R}^{d})}^{2/d} \| u \|_{L_{t}^{2} L_{x}^{\frac{2d}{d - 2}}(J_{k} \times \mathbf{R}^{d})}^{1 + 2/d}.
\end{equation}

\noindent Both $(\ref{2.56})$ and $(\ref{2.57})$ decay to $0$ as $N \nearrow \infty$. $\Box$\vspace{5mm}

\noindent Let

\begin{equation}\label{2.57.1}
C_{0}(N) = \sup((\sup_{J_{k} \subset J} \| u_{|\xi - \xi(a_{k})| \geq N^{1/2}} \|_{S_{\ast}^{0}(J_{k} \times \mathbf{R}^{d})})^{-1/100d}, \frac{1000}{\eta_{1}(d)}),
\end{equation}

\begin{equation}\label{2.57.2}
\eta(N) = \sup(\delta(C_{0}(N))^{-1/100}, C_{0}^{1/100}, 2N^{-1/2}).
\end{equation}

\noindent By lemma $\ref{l2.5}$, $C_{0}(N) \nearrow \infty$, which implies $\eta(N) \searrow 0$. By lemma $\ref{l2.9}$,

\begin{equation}\label{2.58.1}
\| u_{|\xi - \xi(t)| \geq N} \|_{L_{t}^{2} L_{x}^{\frac{2d}{d - 2}}(J \times \mathbf{R}^{d})} \leq C_{2}(m_{0}, d)(\frac{K}{N} + 1)^{1/2} \| u_{|\xi - \xi(t)| \geq \frac{N}{2}} \|_{L_{t}^{\infty} L_{x}^{2}(J \times \mathbf{R}^{d})}
\end{equation}

\begin{equation}\label{2.58.2}
+ C_{2}(m_{0}, d)\sum_{M \leq \eta(N) N} (\frac{M}{N}) \| u_{|\xi - \xi(t)| > M} \|_{L_{t}^{2} L_{x}^{\frac{2d}{d - 2}}(J \times \mathbf{R}^{d})}
\end{equation}

\begin{equation}\label{2.58.3}
+ \delta(C_{0}) \| u_{|\xi - \xi(t)| > \eta(N) N} \|_{L_{t}^{2} L_{x}^{\frac{2d}{d - 2}}(J \times \mathbf{R}^{d})}
\end{equation}

\begin{equation}\label{2.58.4}
+ \left\{
  \begin{array}{ll}
    C_{2}(m_{0}, d) C_{0}^{3/2} (\frac{K}{\eta N})^{1/2} (\sup_{J_{k}} \| u_{|\xi - \xi(a_{k})| \geq \frac{\eta(N) N}{2}} \|_{S_{\ast}^{0}(J_{k} \times \mathbf{R}^{d})}), & \hbox{if $d = 3$;} \\
    C_{2}(m_{0}, d) C_{0}^{4 - 6/d} (\frac{K}{\eta N})^{2/d} \| u_{|\xi - \xi(t)| \geq \eta N} \|_{L_{t}^{2} L_{x}^{\frac{2d}{d - 2}}(J \times \mathbf{R}^{d})}^{1 - 4/d} (\sup_{J_{k}} \| u \|_{S_{\ast}^{0}(J_{k} \times \mathbf{R}^{d})})^{4/d}, & \hbox{if $d \geq 4$.}
  \end{array}
\right.
\end{equation}

\begin{equation}\label{2.58.5}
+ C_{2}(m_{0}, d) (\frac{K}{\eta N})^{1/2} \left\{
                             \begin{array}{ll}
                               \| u_{|\xi - \xi(t)| \geq \eta N} \|_{L_{t}^{\infty} L_{x}^{2}(J \times \mathbf{R}^{3})}, & \hbox{if $d = 3$;} \\
                               \| u_{|\xi - \xi(t)| \geq \eta N} \|_{L_{t}^{\infty} L_{x}^{2}(J \times \mathbf{R}^{3})}^{4/d}, & \hbox{if $d \geq 4$.}
                             \end{array}
                           \right.
\end{equation}

\noindent By theorem $\ref{t2.3}$,

$$(\ref{2.58.2}) \leq 5 C_{2}(m_{0}, d) C_{3}(m_{0}, d) \eta(N)^{1/2} (\frac{K}{N})^{1/2}.$$

$$(\ref{2.58.3}) \leq C_{2}(m_{0}, d) C_{3}(m_{0}, d) \frac{\delta(C_{0}(N))}{\eta(N)^{1/2}} (\frac{K}{N})^{1/2}.$$

\noindent When $d = 3$,

$$(\ref{2.58.4}) \leq C_{2}(m_{0}, d) \frac{C_{0}(N)^{3/2}}{\eta(N)^{1/2}} (\frac{K}{N})^{1/2} (\sup_{J_{k}} \| u_{|\xi - \xi(a_{k})| \geq \frac{\eta(N) N}{2}} \|_{S_{\ast}^{0}(J_{k} \times \mathbf{R}^{d})}),$$

\noindent and when $d \geq 4$,

$$(\ref{2.58.4}) \leq C_{2}(m_{0}, d) C_{3}(m_{0}, d)^{1 - 4/d} \frac{C_{0}(N)^{4 - 6/d}}{\eta(N)^{1/2}} (\frac{K}{N})^{1/2} (\sup_{J_{k}} \| u_{|\xi - \xi(a_{k})| \geq \frac{\eta(N) N}{2}} \|_{S_{\ast}^{0}(J_{k} \times \mathbf{R}^{d})})^{4/d}.$$

\noindent When $d = 3$, let

\begin{equation}\label{2.59}
\aligned
\tilde{\rho}(N) = \frac{C_{2}(m_{0}, d)}{C_{3}(m_{0}, d)} \| u_{|\xi - \xi(t)| \geq \frac{N}{2}} \|_{L_{t}^{\infty} L_{x}^{2}(J \times \mathbf{R}^{d})} 
+ C_{2}(m_{0}, d) (\eta(N))^{1/2} + C_{2}(m_{0}, d) \frac{\delta(C_{0}(N))}{\eta(N)^{1/2}} \\
+ \frac{C_{2}(m_{0}, d)}{C_{3}(m_{0}, d)} \frac{C_{0}(N)^{3/2}}{\eta(N)^{1/2}} (\sup_{J_{k}} \| u_{|\xi - \xi(a_{k})| \geq \frac{\eta(N) N}{2}} \|_{S_{\ast}^{0}(J_{k} \times \mathbf{R}^{d})}),
\endaligned
\end{equation}

\noindent and when $d \geq 4$ let

\begin{equation}\label{2.60}
\aligned
\tilde{\rho}(N) = \frac{C_{2}(m_{0}, d)}{C_{3}(m_{0}, d)} \| u_{|\xi - \xi(t)| \geq \frac{N}{2}} \|_{L_{t}^{\infty} L_{x}^{2}(J \times \mathbf{R}^{d})} 
+ C_{2}(m_{0}, d) (\eta(N))^{1/2} + C_{2}(m_{0}, d) \frac{\delta(C_{0}(N))}{\eta(N)^{1/2}} \\
+ \frac{C_{2}(m_{0}, d)}{C_{3}(m_{0}, d)^{4/d}} \frac{C_{0}(N)^{6 - 4/d}}{\eta(N)^{1/2}} (\sup_{J_{k}} \| u_{|\xi - \xi(a_{k})| \geq \frac{\eta(N) N}{2}} \|_{S_{\ast}^{0}(J_{k} \times \mathbf{R}^{d})})^{4/d}.
\endaligned
\end{equation}

\noindent This implies that for $N \leq K$,

\begin{equation}\label{2.61}
 \| u_{|\xi - \xi(t)| > N} \|_{L_{t}^{2} L_{x}^{\frac{2d}{d - 2}}(J \times \mathbf{R}^{d})} \leq C_{3}(m_{0}, d) \tilde{\rho}(N) (\frac{K}{N})^{1/2}.
\end{equation}

\noindent Lemma $\ref{l2.5}$, $(\ref{2.57.1})$, $(\ref{2.57.2})$ imply $\tilde{\rho}(N) \rightarrow 0$ as $N \rightarrow \infty$. Taking $\rho(N) = \inf(1, \tilde{\rho}(N))$ proves the theorem. $\Box$\vspace{5mm}

\noindent \textbf{Remark:} These estimates also hold for $u$ a minimal mass blowup solution to the focusing initial value problem $(\ref{0.2.4.1})$.


\section{$\int_{0}^{\infty} N(t)^{3} dt = \infty$} 

\noindent We will defeat this scenario by proving a frequency localized interaction Morawetz estimate. The interaction Morawetz estimate was proved for solutions to the defocusing nonlinear Schr{\"o}dinger equation in \cite{CKSTT2} when $d = 3$, and in \cite{TVZ} for dimensions $d \geq 4$. The interaction Morawetz estimate was proved by taking the tensor product of two solutions to $(\ref{0.1})$. Let $x$ refer to the first $d$ variables in $\mathbf{R}^{d} \times \mathbf{R}^{d}$ and $y$ refer to the second $d$ variables. We adopt the convention of summing over repeated indices. Let $M(t)$ be the Morawetz action

\begin{equation}\label{4.0.0.0}
M(t) = \int_{\mathbf{R}^{d} \times \mathbf{R}^{d}} \frac{(x - y)_{j}}{|x - y|} Im[ \bar{u}(t,x) \bar{u}(t,y) \partial_{j}(u(t,x) u(t,y))] dx dy.
\end{equation}

\noindent \cite{CKSTT2} proved

\begin{equation}\label{4.0.0}
\| u \|_{L_{t,x}^{4}(I \times \mathbf{R}^{3})}^{4} \lesssim \int_{I} \partial_{t} M(t) dt \lesssim \sup_{t \in I} |M(t)| \lesssim \| u \|_{L_{t}^{\infty} L_{x}^{2}(I \times \mathbf{R}^{3})}^{3} \| u \|_{L_{t}^{\infty} \dot{H}_{x}^{1}(I \times \mathbf{R}^{3})}.
\end{equation}

\noindent \cite{TVZ} proved

\begin{equation}\label{4.0.0.1}
\aligned
\int_{I} \int_{\mathbf{R}^{d} \times \mathbf{R}^{d}} (-\Delta \Delta |x - y|) |u(t,x)|^{2} |u(t,y)|^{2} dx dy dt \lesssim \int_{I} \partial_{t} M(t) dt \\ \lesssim \sup_{t \in I} |M(t)| \lesssim \| u \|_{L_{t}^{\infty} L_{x}^{2}(I \times \mathbf{R}^{d})}^{3} \| u \|_{L_{t}^{\infty} \dot{H}_{x}^{1}(I \times \mathbf{R}^{d})}.
\endaligned
\end{equation}

\noindent Additionally, the quantities $\| u \|_{L_{t,x}^{4}(I \times \mathbf{R}^{3})}$ and $$\int_{I} \int_{\mathbf{R}^{d} \times \mathbf{R}^{d}} (-\Delta \Delta |x - y|) |u(t,x)|^{2} |u(t,y)|^{2} dx dy dt$$ are invariant under the transformation $u \mapsto e^{ix \cdot \xi(t)} u$. We will show that $M(t)$ is also Galilean invariant. See \cite{PV} for more information.\vspace{5mm}

\noindent Indeed, let

\begin{equation}\label{4.0.0.2}
\tilde{M}(t) = \int_{\mathbf{R}^{d} \times \mathbf{R}^{d}} \frac{(x - y)_{j}}{|x - y|} Im[\overline{u}(t,x) \overline{u}(t,y) (\partial_{j} - i \xi_{j}(t)) u(t,x) u(t,y)] dx dy.
\end{equation}

\noindent Then $$\int_{\mathbf{R}^{d} \times \mathbf{R}^{d}} \frac{(x - y)_{j}}{|x - y|} Im[\overline{u}(t,x) \overline{u}(t,y) (i \xi_{j}(t)) u(t,x) u(t,y)] dx dy$$

$$ = \int_{\mathbf{R}^{d} \times \mathbf{R}^{d}} \xi_{j}(t) \frac{(x - y)_{j}}{|x - y|} |u(t,x)|^{2} |u(t,y)|^{2} dx dy.$$

\noindent Because $|u(t,x)|^{2} |u(t,y)|^{2}$ is even in $x - y$ and $\frac{(x - y)_{j}}{|x - y|}$ is odd in $x - y$, $M(t) = \tilde{M}(t)$.\vspace{5mm}

\noindent We will not use these estimates directly, instead, we use a frequency localized interaction Morawetz estimate. \cite{CKSTT4} introduced a frequency localized version of $(\ref{4.0.0})$ for the energy critical nonlinear Schr{\"o}dinger equation on $\mathbf{R}^{3}$ to prove global well-posedness and scattering. In that case $u(t) \in \dot{H}^{1}(\mathbf{R}^{3})$, so the Morawetz estimates were localized to high frequencies. Here $u(t) \in L^{2}(\mathbf{R}^{3})$, so we localize to low frequencies. In the energy critical case, $d = 3$, the $L_{t,x}^{4}$ norm scales like $$\int_{I} N(t)^{-1} dt,$$ while in the mass critical case the $L_{t,x}^{4}$ norm scales like $$\int_{I} N(t)^{3} dt.$$

\noindent This method also has a great deal in common with the almost Morawetz estimates frequently used in conjunction with the I-method. (See \cite{CGT}, \cite{CR}, and \cite{D} for the two dimensional case, and \cite{D1} in the three dimensional case.)\vspace{5mm}

\noindent Let C be a fixed constant and let $m(\xi)$ be the smooth, radial Fourier multiplier,

\begin{equation}\label{4.0}
m(\xi) = \left\{
           \begin{array}{ll}
             1, & \hbox{$|\xi| \leq CK$;} \\
             0, & \hbox{$|\xi| > 2CK$.}
           \end{array}
         \right.
\end{equation}

\begin{theorem}\label{t4.1}
Suppose $J$ is a compact interval with $N(t) \leq 1$ and $\int_{J} N(t)^{3} dt = K$. Then if $u$ is a minimal mass blowup solution to $(\ref{0.1})$, $\mu = +1$,

\begin{equation}\label{4.1}
\int_{J} \int_{\mathbf{R}^{d} \times \mathbf{R}^{d}} (-\Delta \Delta |x - y|) |\I u(t,x)|^{2} |\I u(t,y)|^{2} dx dy dt \lesssim_{m_{0}, d} o(K).
\end{equation}

\noindent $o(K)$ is a quantity with $\lim_{K \rightarrow \infty} \frac{o(K)}{K} = 0$.
\end{theorem}

\noindent \textbf{Remark:} The interaction Morawetz estimates of \cite{TVZ}, \cite{PV}, \cite{CGT1}, and \cite{CKSTT2} rely heavily on $\mu = +1$. When $\mu = -1$ the interaction Morawetz estimates are no longer positive definite, and therefore do not give an estimate of the form $(\ref{4.0.0})$. This is the main obstacle to extending our methods from the defocusing case to the focusing case.\vspace{5mm}

\noindent \textbf{Remark:} Since $J$ is a compact interval and $N(t) \leq 1$, $$\| u \|_{L_{t,x}^{\frac{2(d + 2)}{d}}(J \times \mathbf{R}^{d})} < \infty.$$ This means $J$ can be partitioned into a finite number of intervals $J_{k}$ with $\| u \|_{L_{t,x}^{\frac{2(d + 2)}{d}}(J_{k} \times \mathbf{R}^{d})} = \epsilon$. By lemma $\ref{l1.4}$, $$\sum_{J_{k}} N(J_{k}) \sim \int_{J} N(t)^{3} dt.$$ Therefore theorem $\ref{t4.1}$ is good enough to exclude the scenario $\int_{0}^{\infty} N(t)^{3} dt = \infty$.\vspace{5mm}

\noindent \textbf{Remark:} For the rest of this section we will simply write $\lesssim$ and understand that this refers to $\lesssim_{m_{0}, d}$.\vspace{5mm}

\begin{theorem}\label{t4.2}
If theorem $\ref{t4.1}$ is true, then there does not exist a minimal mass blowup solution to $(\ref{0.1})$ with $N(t) \leq 1$, $\mu = +1$, and $$\int_{0}^{\infty} N(t)^{3} dt = \infty.$$
\end{theorem}

\noindent \emph{Proof of Theorem $\ref{t4.2}$:} We want $\| \I u(t) \|_{L^{2}(\mathbf{R}^{d})}$ to be very close to $\| u(t) \|_{L^{2}(\mathbf{R}^{d})}$ for all $t$. Therefore, make a Galilean transformation so that $\xi(0) = 0$. Consider $d = 3$ and $d \geq 4$ separately.\vspace{5mm}

\noindent \textbf{Case 1, $d = 3$:} In this case we need a local well-posedness result.

\begin{lemma}\label{l4.2}
Suppose $J_{1}$ is an interval with $\| \I u \|_{L_{t,x}^{10/3}(J_{1} \times \mathbf{R}^{3})} = \frac{\epsilon}{2}$, $C$ is very large, and $|\xi(t)| \lesssim_{d} K$. Then $\| u \|_{L_{t,x}^{10/3}(J_{1} \times \mathbf{R}^{3})} \leq \frac{3 \epsilon}{4}$.
\end{lemma}

\noindent \emph{Proof:} Without loss of generality let $J_{1} = [0, T]$. By Duhamel's formula and Strichartz estimates,

\begin{equation}\label{4.4}
\aligned
\| u \|_{S^{0}(J_{1} \times \mathbf{R}^{3})} &\lesssim \| u_{0} \|_{L^{2}(\mathbf{R}^{3})} + \| \I u \|_{L_{t,x}^{10/3}(J_{1} \times \mathbf{R}^{3})}^{7/3} + \\ &\| (1 - \I)u \|_{L_{t}^{\infty} L_{x}^{2}(J_{1} \times \mathbf{R}^{3})}^{4/3} \| (1 - \I) u \|_{L_{t}^{2} L_{x}^{6}(J_{1} \times \mathbf{R}^{3})} \\
&\lesssim \| u_{0} \|_{L^{2}(\mathbf{R}^{3})} + \epsilon^{4/d} + \| (1 - \I) u \|_{L_{t}^{\infty} L_{x}^{2}(J_{1} \times \mathbf{R}^{3})}^{4/3} \| u \|_{S^{0}(J_{1} \times \mathbf{R}^{3})}.
\endaligned
\end{equation}

\noindent Since $\| u_{0} \|_{L^{2}(\mathbf{R}^{3})} \lesssim 1$, $\| (1 - \I) u \|_{L_{t}^{\infty} L_{x}^{2}(J_{1} \times \mathbf{R}^{3})}$ sufficiently small implies $\| u \|_{S^{0}(J \times \mathbf{R}^{3})} \lesssim 1$ by continuity. Interpolating $\| (1 - \I) u \|_{L_{t}^{2} L_{x}^{6}(J_{1} \times \mathbf{R}^{3})} \lesssim 1$ with $\| (1 - \I)u \|_{L_{t}^{\infty} L_{x}^{2}(J_{1} \times \mathbf{R}^{3})} \leq \delta(\epsilon)$ for $\delta(\epsilon) > 0$ sufficiently small implies $\| u \|_{L_{t,x}^{10/3}(J_{1} \times \mathbf{R}^{3})} \leq \frac{3 \epsilon}{4}$. By $(\ref{6.3})$, $|\xi(t)| \lesssim_{d} K$, so we can choose $C(\delta, d)$ sufficiently large so that $$\| u_{> \frac{CK}{2}} \|_{L_{t}^{\infty} L_{x}^{2}(J_{1} \times \mathbf{R}^{3})} \leq \delta(\epsilon).$$ $\Box$\vspace{5mm}

\noindent \textbf{Remark:} By lemma $\ref{l1.2}$, if $\int_{J} N(t)^{3} dt = K$, then for any $t_{1}, t_{2} \in J$, $|\xi(t_{1}) - \xi(t_{2})| \lesssim_{d} K$. Therefore if $\int_{J} N(t)^{3} dt = K$ we can make a Galilean transformation so that $|\xi(t)| \lesssim_{d} K$ on $J$.\vspace{5mm}

\noindent Now take a subinterval $J_{k}$ with $\| u \|_{L_{t,x}^{10/3}(J_{k} \times \mathbf{R}^{3})} = \epsilon$. Lemma $\ref{l4.2}$ implies that $\| \I u \|_{L_{t,x}^{10/3}(J_{k} \times \mathbf{R}^{3})} \geq \frac{\epsilon}{2}$. From $(\ref{1.12.1})$,

\begin{equation}\label{4.4.1}
\int_{J_{k}} N(t)^{2} dt \lesssim \int_{J_{k}} \int_{\mathbf{R}^{3}} |u(t,x)|^{10/3} dx dt \lesssim \epsilon^{10/3}.
\end{equation}

\noindent By lemma $\ref{l1.1}$, $N(t_{1}) \sim N(t_{2})$ on $J_{k}$, so

$$|J_{k}| \lesssim \frac{\epsilon^{10/3}}{N(J_{k})^{2}}.$$ By Holder's inequality,

\begin{equation}\label{4.5}
\| \I u \|_{L_{t}^{8/3} L_{x}^{4}(J_{k} \times \mathbf{R}^{3})} \lesssim (\frac{1}{N(J_{k})^{2}})^{1/8} \| \I u \|_{L_{t,x}^{4}(J_{k} \times \mathbf{R}^{d})}.
\end{equation}

\noindent This implies

\begin{equation}\label{4.6}
N(J_{k}) \| \I u \|_{L_{t}^{8/3} L_{x}^{4}(J_{k} \times \mathbf{R}^{3})}^{4} \lesssim \| \I u \|_{L_{t,x}^{4}(J_{k} \times \mathbf{R}^{3})}^{4}.
\end{equation}

\noindent By interpolation if $\| \I u \|_{L_{t,x}^{10/3}(J \times \mathbf{R}^{3})} \geq \frac{\epsilon}{2}$ and $\| \I u \|_{L_{t}^{\infty} L_{x}^{2}(J_{k} \times \mathbf{R}^{3})} \lesssim 1$, then $\| \I u \|_{L_{t}^{8/3} L_{x}^{4}(J_{k} \times \mathbf{R}^{3})} \gtrsim \epsilon^{5/4}$, so

\begin{equation}\label{4.7}
\aligned
\int_{J} N(t)^{3} dt \sim \sum_{J_{k}} N(J_{k}) \lesssim \sum_{J_{k}} N(J_{k}) \| \I u \|_{L_{t,x}^{8/3}(J_{k} \times \mathbf{R}^{3})}^{4} \\ \lesssim \sum_{J_{k}} \| \I u \|_{L_{t,x}^{4}(J_{k} \times \mathbf{R}^{3})}^{4} = \int_{J} \int_{\mathbf{R}^{3}} |\I u(t,x)|^{4} dx dt.
\endaligned
\end{equation}

\noindent When $d = 3$, $$(-\Delta \Delta |x - y|) = 4 \pi \delta(|x - y|).$$ Therefore $$\int_{\mathbf{R}^{3} \times \mathbf{R}^{3}} (-\Delta \Delta |x - y|) |\I u(t,y)|^{2} |\I u(t,x)|^{2} dx dy = \int_{\mathbf{R}^{3}} |\I u(t,x)|^{4} dx.$$ Now if $$\int_{0}^{T} N(t)^{3} dt = K,$$ then by theorem $\ref{t4.1}$,

\begin{equation}\label{4.8}
K \lesssim_{d} \int_{0}^{T} \int_{\mathbf{R}^{3}} |Iu(t,x)|^{4} dx dt \lesssim_{d} o(K).
\end{equation}

\noindent This gives a contradiction if $K$ is sufficiently large. When $\int_{0}^{\infty} N(t)^{3} dt = \infty$ we can always a suitable $T$.\vspace{5mm}

\noindent \textbf{Case 2, $d \geq 4$:}

\begin{equation}\label{4.9}
\aligned
&\int_{J} \int_{\mathbf{R}^{d} \times \mathbf{R}^{d}} (-\Delta \Delta |x - y|) |\I u(t,x)|^{2} |\I u(t,y)|^{2} dx dy dt \\ =  &\int_{J} \int_{\mathbf{R}^{d} \times \mathbf{R}^{d}} (\frac{4(d - 1)(d - 3)}{|x - y|^{3}}) |\I u(t,x)|^{2} |\I u(t,y)|^{2} dx dy dt.
\endaligned
\end{equation}

\noindent Let $\eta = \frac{m_{0}^{2}}{1000}$.

\begin{equation}\label{4.10}
\int_{|x - x(t)| \leq \frac{C(\eta)}{N(t)}} |u(t,x)|^{2} dx \geq m_{0}^{2} - \eta.
\end{equation}

\noindent Also,

\begin{equation}\label{4.11}
\int_{|\xi - \xi(t)| > \frac{C(\eta) N(t)}{2}} |\hat{u}(t,\xi)|^{2} d\xi \leq \eta.
\end{equation}

\noindent Therefore, for $K \geq 1$, $C$ sufficiently large,

\begin{equation}\label{4.12}
\int_{|x - x(t)| \leq \frac{C(\eta)}{N(t)}} |\I u(t,x)|^{2} dx \geq \frac{m_{0}^{2}}{2}.
\end{equation}

\noindent Of course, for the same $x(t) \in \mathbf{R}^{d}$ we also have

\begin{equation}\label{4.12.1}
\int_{|y - x(t)| \leq \frac{C(\eta)}{N(t)}} |\I u(t,y)|^{2} dx \geq \frac{m_{0}^{2}}{2}.
\end{equation}

\noindent Therefore, because $N(t) \leq 1$,

$$N(t)^{3} \lesssim N(t)^{3} (\int_{|x - x(t)| \leq \frac{C(\eta)}{N(t)}} |\I u(t,x)|^{2} dx)(\int_{|y - x(t)| \leq \frac{C(\eta)}{N(t)}} |\I u(t,y)|^{2} dy)$$

$$\lesssim N(t)^{3} \int_{|x - y| \leq \frac{2 C(\eta)}{N(t)}} |\I u(t,x)|^{2} |\I u(t,y)|^{2} dx dy$$

$$\lesssim \int_{\mathbf{R}^{d} \times \mathbf{R}^{d}} \frac{1}{|x - y|^{3}} |\I u(t,x)|^{2} |\I u(t,y)|^{2} dx dy.$$

\noindent Once again, this implies that for a compact interval $J$,

\begin{equation}\label{4.13}
K = \int_{J} N(t)^{3} dt \lesssim_{d} \int_{J} \int_{\mathbf{R}^{d} \times \mathbf{R}^{d}} (\frac{1}{|x - y|})^{3} |\I u(t,x)|^{2} |\I u(t,y)|^{2} dx dy dt \lesssim o(K).
\end{equation}

\noindent This gives a contradiction for $K$ sufficiently large. $\Box$\vspace{5mm}

\noindent All that is left to do is to prove theorem $\ref{t4.1}$, which will occupy the remainder of the section. We begin by estimating the error for the truncated Morawetz estimates. For the rest of the section $C(\epsilon, m_{0}, d)$ will be a fixed constant so that $(\ref{4.12})$ is satisfied, $|\xi(t)| \leq \frac{CK}{1000}$ on $J$ if $\int_{J} N(t)^{3} dt = K$, and $\| \I u \|_{L_{t,x}^{10/3}(J_{1} \times \mathbf{R}^{3})} \leq \frac{\epsilon}{2}$ implies $\| u \|_{L_{t,x}^{10/3}(J_{1} \times \mathbf{R}^{3})} \leq \frac{3 \epsilon}{4}$.\vspace{5mm}

\begin{theorem}\label{t4.3}
Let $a(x,y) = |x - y|$. Define the interaction Morawetz quantity

\begin{equation}\label{4.14}
M(t) = \int a_{j}(x,y) Im[\overline{\I u(t,x)} \overline{\I u(t,y)} (\partial_{j} - i \xi_{j}(t)) (\I u(t,x) \I u(t,y))] dx dy.
\end{equation}

\noindent Then for $\mu = +1$,

\begin{equation}\label{4.15}
\int_{0}^{T} \int_{\mathbf{R}^{d} \times \mathbf{R}^{d}} (-\Delta \Delta |x - y|) |\I u(t,x)|^{2} |\I u(t,y)|^{2} dx dy dt \lesssim  o(K).
\end{equation}

\end{theorem}

\noindent \textbf{Remark:} We adopt the usual convention of summing over repeated indices.\vspace{5mm}

\noindent \emph{Proof:} First take $M(t)$.

$$\int_{\mathbf{R}^{d} \times \mathbf{R}^{d}} \frac{(x - y)_{j}}{|x - y|} Im[\overline{\I u}(t,x) \overline{\I u}(t,y) (\partial_{j} - i \xi_{j}(t)) \I u(t,x) \I u(t,y)] dx dy$$

$$\lesssim \| \I u \|_{L_{t}^{\infty} L_{x}^{2}([0, T] \times \mathbf{R}^{d})}^{3} \| (\nabla - i \xi(t)) \I u \|_{L_{t}^{\infty} L_{x}^{2}([0, T] \times \mathbf{R}^{d})} \lesssim o(K).$$

\noindent We estimate $\| \I u \|_{L_{t}^{\infty} L_{x}^{2}([0, T] \times \mathbf{R}^{d})}$ by conservation of mass and $\| (\nabla - i \xi(t)) \I u \|_{L_{t}^{\infty} L_{x}^{2}([0, T] \times \mathbf{R}^{d})}$ by $(\ref{6.3})$ and $N(t) \leq 1$.\vspace{5mm}

\noindent Since $\I$ is a Fourier multiplier,

\begin{equation}\label{4.15.1}
\partial_{t} (\I u) = i \Delta \I u - i | \I u|^{4/d} (\I u) + i |\I u|^{4/d} (\I u) - i \I(|u|^{4/d} u).
\end{equation}

\noindent If we had only $$\partial_{t} (\I u) = i \Delta (\I u) - i |\I u|^{4/d} (\I u)$$ then the proof of theorem $\ref{t4.3}$ would be complete. We could copy the arguments from \cite{CKSTT2} and \cite{TVZ} exactly, replacing $u$ with $\I u$. Instead, it is necessary to deal with the error terms that arise from the fact that $|\I u|^{4/d} (\I u) - \I(|u|^{4/d} u) \neq 0$, and prove these error terms are $\lesssim o(K)$. Let $x$ denote the first $d$ variables in $\mathbf{R}^{d} \times \mathbf{R}^{d}$ and $y$ the second $d$ variables. We have the error

\begin{equation}\label{4.15.2}
\aligned
\mathcal E = \int_{0}^{T} \int_{\mathbf{R}^{d} \times \mathbf{R}^{d}} a_{j}(x,y) |\I u(t,y)|^{2} \\ Re \{ [\I(|u|^{4/d} \bar{u})(t,x) - |\I u|^{4/d} (\overline{\I u})(t,x)]  (\partial_{j} - i \xi_{j}(t)) \I u(t,x) \} dx dy dt
\endaligned
\end{equation}

\begin{equation}\label{4.15.3}
\aligned
+ \int_{0}^{T} \int_{\mathbf{R}^{d} \times \mathbf{R}^{d}} a_{j}(x,y) |\I u(t,y)|^{2} \\ Re \{ \overline{\I u}(t,x) (\partial_{j} - i \xi_{j}(t)) [|\I u|^{4/d} (\I u)(t,x) - \I(|u|^{4/d} u)(t,x)] \} dx dy dt
\endaligned
\end{equation}

\begin{equation}\label{4.15.4}
\aligned
+ \int_{0}^{T} \int_{\mathbf{R}^{d} \times \mathbf{R}^{d}} a_{j}(x,y) Re[&[\overline{\I u}(t,x) (\partial_{j} - i \xi_{j}(t)) \I u(t,x)] \\ &[\I(|u|^{4/d} \bar{u})(t,y) \I u(t,y)  - \I(|u|^{4/d} u)(t,y) \overline{\I u}(t,y)]] dx dy dt.
\endaligned
\end{equation}

\noindent Now we need some intermediate lemmas.

\begin{lemma}\label{l4.5.1}
Suppose $u$ satisfies

\begin{equation}\label{4.16.1}
\| P_{|\xi - \xi(t)| > N} u \|_{L_{t}^{2} L_{x}^{\frac{2d}{d - 2}}([0, T] \times \mathbf{R}^{d})} \lesssim_{m_{0}, d} \rho(N) ((\frac{K}{N})^{1/2} + 1),
\end{equation}

\noindent $\rho(N) \leq 1$, $\rho(N) \rightarrow 0$ as $N \rightarrow \infty$, $|\xi(t)| \leq \frac{CK}{1000}$. Then for any $1/2 < s \leq 1$,

\begin{equation}\label{4.16.2}
\| |\nabla|^{s} e^{-ix \cdot \xi(t)} \I u \|_{L_{t}^{2} L_{x}^{\frac{2d}{d - 2}}([0, T] \times \mathbf{R}^{d})} \lesssim o(K^{s}).
\end{equation}

\end{lemma}

\noindent \emph{Proof:} $$\| |\nabla|^{s} (e^{-ix \cdot \xi(t)} \I u) \|_{L_{t}^{2} L_{x}^{\frac{2d}{d - 2}}([0, T] \times \mathbf{R}^{d})} \lesssim \sum_{N \leq 2CK} N^{s} \| P_{N}(e^{-ix \cdot \xi(t)} \I u) \|_{L_{t}^{2} L_{x}^{\frac{2d}{d - 2}}([0, T] \times \mathbf{R}^{d})}.$$

$$\lesssim \sum_{N \leq 2CK} N^{s} \| P_{N}(e^{-ix \cdot \xi(t)} u) \|_{L_{t}^{2} L_{x}^{\frac{2d}{d - 2}}([0, T] \times \mathbf{R}^{d})} \lesssim \sum_{N \leq 2CK} N^{s} \rho(N) (\frac{K}{N})^{1/2} \lesssim o(K^{s}).$$

 $\Box$

\begin{lemma}\label{l4.5.2}
Suppose $u$ satisfies the hypotheses of lemma $\ref{l4.5.1}$. Then for $1/2 < s \leq 1$,

\begin{equation}\label{4.16.3}
\| |\nabla|^{s} (e^{-ix \cdot \xi(t)} \I(|u|^{4/d} u)) \|_{L_{t}^{2} L_{x}^{\frac{2d}{d + 2}}([0, T] \times \mathbf{R}^{d})} \lesssim_{m_{0}, d} K^{s}.
\end{equation}

\end{lemma}

\noindent \emph{Proof:} Again make a Littlewood-Paley decomposition.

\begin{equation}\label{4.16.4}
\aligned
\| |\nabla|^{s} (e^{-ix \cdot \xi(t)} \I(|u|^{4/d} u)) \|_{L_{t}^{2} L_{x}^{\frac{2d}{d + 2}}([0, T] \times \mathbf{R}^{d})} \\ \leq \sum_{N \leq 2CK} N^{s} \| P_{N}(e^{-ix \cdot \xi(t)} \I(|u|^{4/d} u)) \|_{L_{t}^{2} L_{x}^{\frac{2d}{d + 2}}([0, T] \times \mathbf{R}^{d})}
\endaligned
\end{equation}

$$\sum_{N \leq \frac{CK}{4}} N^{s} \| P_{N} (e^{-ix \cdot \xi(t)} \I(|u|^{4/d} u)) \|_{L_{t}^{2} L_{x}^{\frac{2d}{d + 2}}([0, T] \times \mathbf{R}^{d})}$$ $$= \sum_{N \leq \frac{CK}{4}} N^{s} \| P_{N} (e^{-ix \cdot \xi(t)} (|u|^{4/d} u)) \|_{L_{t}^{2} L_{x}^{\frac{2d}{d + 2}}([0, T] \times \mathbf{R}^{d})}.$$

\noindent By Bernstein's inequality,

$$\| P_{N}(|P_{\leq N} (e^{-ix \cdot \xi(t)} u)|^{4/d} (P_{\leq N} (e^{-ix \cdot \xi(t)} u))) \|_{L_{t}^{2} L_{x}^{\frac{2d}{d + 2}}([0, T] \times \mathbf{R}^{d})}$$ $$\lesssim \frac{1}{N} \| u \|_{L_{t}^{\infty} L_{x}^{2}([0, T] \times \mathbf{R}^{d})}^{4/d} \| \nabla (P_{\leq N} (e^{-ix \cdot \xi(t)} u))) \|_{L_{t}^{2} L_{x}^{\frac{2d}{d - 2}}([0, T] \times \mathbf{R}^{d})} \lesssim \frac{K^{1/2}}{N^{1/2}}.$$

\noindent By Holder's inequality, conservation of mass,

$$\| |P_{> N} (e^{-ix \cdot \xi(t)} u)| |u|^{4/d} \|_{L_{t}^{2} L_{x}^{\frac{2d}{d + 2}}([0, T] \times \mathbf{R}^{d})} \lesssim \frac{K^{1/2}}{N^{1/2}}.$$

\noindent Therefore, for $N \leq \frac{CK}{4}$,

$$\| P_{N} (e^{-ix \cdot \xi(t)} \I(|u|^{4/d} u)) \|_{L_{t}^{2} L_{x}^{\frac{2d}{d + 2}}([0, T] \times \mathbf{R}^{d})} \lesssim \frac{K^{1/2}}{N^{1/2}}.$$

\noindent Meanwhile,

$$\| P_{\geq \frac{CK}{4}} (e^{-ix \cdot \xi(t)} \I(|u|^{4/d} u)) \|_{L_{t}^{2} L_{x}^{\frac{2d}{d + 2}}([0, T] \times \mathbf{R}^{d})}$$

$$\leq \| P_{\geq \frac{CK}{5}}  (|u|^{4/d} u) \|_{L_{t}^{2} L_{x}^{\frac{2d}{d + 2}}([0, T] \times \mathbf{R}^{d})}$$

$$\leq \| P_{\geq \frac{CK}{8}}  (|e^{-ix \cdot \xi(t)} u|^{4/d} (e^{-ix \cdot \xi(t)} u)) \|_{L_{t}^{2} L_{x}^{\frac{2d}{d + 2}}([0, T] \times \mathbf{R}^{d})}.$$

\noindent Again combining Bernstein's inequality, conservation of mass, and Holder's inequality,

$$\| P_{\geq \frac{CK}{8}}  (|e^{-ix \cdot \xi(t)} u|^{4/d} (e^{-ix \cdot \xi(t)} u)) \|_{L_{t}^{2} L_{x}^{\frac{2d}{d + 2}}([0, T] \times \mathbf{R}^{d})} \lesssim 1.$$

\noindent Therefore $(\ref{4.16.4}) \lesssim K^{s}$. $\Box$\vspace{5mm}

\begin{lemma}\label{l4.5.3}
Suppose $u$ satisfies the hypotheses of lemma $\ref{l4.5.1}$. Then

\begin{equation}\label{4.16.5}
\| \I(|u|^{4/d} u) - |\I u|^{4/d} (\I u) \|_{L_{t}^{2} L_{x}^{\frac{2d}{d + 2}}([0, T] \times \mathbf{R}^{d})} \lesssim_{m_{0}, d} 1.
\end{equation}

\end{lemma}

\noindent \emph{Proof:} By lemma $\ref{l4.5.2}$, $$\| \nabla e^{-ix \cdot \xi(t)} \I(|u|^{4/d} u) \|_{L_{t}^{2} L_{x}^{\frac{2d}{d + 2}}([0, T] \times \mathbf{R}^{d})} \lesssim K.$$ Also, by the chain rule and conservation of mass,

$$\| \nabla e^{-ix \cdot \xi(t)} (|\I u|^{4/d} (\I u)) \|_{L_{t}^{2} L_{x}^{\frac{2d}{d + 2}}([0, T] \times \mathbf{R}^{d})} \lesssim_{m_{0}, d} \| \nabla e^{-ix \cdot \xi(t)} \I u \|_{L_{t}^{2} L_{x}^{\frac{2d}{d - 2}}([0, T] \times \mathbf{R}^{d})}.$$

\noindent Since $\I = 1$ on $|\xi| \leq CK$,

$$\| \nabla e^{-ix \cdot \xi(t)} \I(P_{\leq \frac{CK}{4}} u) \|_{L_{t}^{2} L_{x}^{\frac{2d}{d - 2}}([0, T] \times \mathbf{R}^{d})} = \| \nabla \I e^{-ix \cdot \xi(t)} (P_{\leq \frac{CK}{4}} u) \|_{L_{t}^{2} L_{x}^{\frac{2d}{d - 2}}([0, T] \times \mathbf{R}^{d})} \lesssim K.
$$

\noindent The last inequality follows from lemma $\ref{l4.5.1}$.

$$\| \nabla e^{-ix \cdot \xi(t)} \I(P_{\geq \frac{CK}{4}} u) \|_{L_{t}^{2} L_{x}^{\frac{2d}{d - 2}}([0, T] \times \mathbf{R}^{d})}
\lesssim K \| u_{|\xi - \xi(t)| > \frac{CK}{8}} \|_{L_{t}^{2} L_{x}^{\frac{2d}{d - 2}}([0, T] \times \mathbf{R}^{d})} \lesssim K.$$

\noindent Therefore, by Bernstein's inequality,

\begin{equation}\label{4.16.6}
\| P_{> \frac{CK}{4}} e^{-ix \cdot \xi(t)} [\I(|u|^{4/d} u) - |\I u|^{4/d} (\I u)] \|_{L_{t}^{2} L_{x}^{\frac{2d}{d + 2}}([0, T] \times \mathbf{R}^{d})} \lesssim 1.
\end{equation}

\noindent On the other hand, by $|\xi(t)| \leq \frac{CK}{1000}$ and Holder's inequality,

$$
\| P_{\leq \frac{CK}{4}} e^{-ix \cdot \xi(t)} [\I(|u|^{4/d} u) - |\I u|^{4/d} (\I u)] \|_{L_{t}^{2} L_{x}^{\frac{2d}{d + 2}}([0, T] \times \mathbf{R}^{d})}$$ $$\leq \| |u|^{4/d} u - |\I u|^{4/d} (\I u) \|_{L_{t}^{2} L_{x}^{\frac{2d}{d + 2}}([0, T] \times \mathbf{R}^{d})}$$

$$\lesssim \| |P_{> \frac{CK}{4}} u | |u|^{4/d} \|_{L_{t}^{2} L_{x}^{\frac{2d}{d + 2}}([0, T] \times \mathbf{R}^{d})} \lesssim \| u_{|\xi - \xi(t)| \geq \frac{CK}{8}} \|_{L_{t}^{2} L_{x}^{\frac{2d}{d - 2}}([0, T] \times \mathbf{R}^{d})} \lesssim 1.$$

\noindent Therefore the proof is complete. $\Box$\vspace{5mm}

\noindent We are now ready to estimate the first term in $\mathcal E$.

\begin{corollary}\label{c4.5.4}
\begin{equation}\label{4.16.7}
 (\ref{4.15.2}) \lesssim o(K).
\end{equation}

\end{corollary}

\noindent \emph{Proof:} Because $\frac{(x - y)_{j}}{|x - y|}$ is uniformly bounded on $\mathbf{R}^{d} \times \mathbf{R}^{d}$, by lemmas $\ref{l4.5.1}$, $\ref{l4.5.3}$,

$$(\ref{4.15.2}) \lesssim \| \I u \|_{L_{t}^{\infty} L_{x}^{2}([0, T] \times \mathbf{R}^{d})}^{2} \| e^{ix \cdot \xi(t)} \nabla (e^{-ix \cdot \xi(t)} \I u) \|_{L_{t}^{2} L_{x}^{\frac{2d}{d - 2}}([0, T] \times \mathbf{R}^{d})}$$ $$\times \| \I(|u|^{4/d} u) - |\I u|^{4/d} (\I u) \|_{L_{t}^{2} L_{x}^{\frac{2d}{d + 2}}([0, T] \times \mathbf{R}^{d})} \lesssim o(K).$$ $\Box$\vspace{5mm}

\noindent In order to estimate $(\ref{4.15.3})$ and $(\ref{4.15.4})$ we need one additional lemma.

\begin{lemma}\label{l4.5.5}
Suppose $K(x)$ is a kernel,

\begin{equation}\label{4.17.1}
|K(x)| \lesssim_{d} 1,
\end{equation}

\noindent and

\begin{equation}\label{4.17.2}
|\nabla K(x)| \lesssim_{d} \frac{1}{|x|}.
\end{equation}

\noindent Let

\begin{equation}\label{4.17.3}
F(x) = \int K(x - y) \cdot (\nabla f(y)) g(y) dy.
\end{equation}

\noindent Then $F(x) = G(x) + H(x)$, where for $\frac{1}{p} + \frac{1}{p'} = 1$,

\begin{equation}\label{4.17.4}
\| G \|_{L_{x}^{\infty}(\mathbf{R}^{d})} \lesssim_{d} \| \nabla g \|_{L_{x}^{p}(\mathbf{R}^{d})} \| f \|_{L_{x}^{p'}(\mathbf{R}^{d})},
\end{equation}

\begin{equation}\label{4.17.5}
\| H \|_{L_{x}^{\frac{6d}{5}}(\mathbf{R}^{d})} \lesssim_{d} \| |\nabla|^{2/3} g \|_{L_{x}^{\frac{2d}{d - 1}}(\mathbf{R}^{d})} \| f \|_{L_{x}^{\frac{2d}{d + 2}}(\mathbf{R}^{d})},
\end{equation}

\noindent and

\begin{equation}\label{4.17.6}
\| H \|_{L_{x}^{3d}(\mathbf{R}^{d})} \lesssim_{d} \| |\nabla|^{2/3} g \|_{L_{x}^{\frac{2d}{d + 2}}(\mathbf{R}^{d})} \| f \|_{L_{x}^{\frac{2d}{d - 2}}(\mathbf{R}^{d})}.
\end{equation}
\end{lemma}

\noindent \emph{Proof:} This is proved by integration by parts and the Hardy-Littlewood-Sobolev inequality.

$$\int K(x - y) \cdot (\nabla f(y)) g(y) dy = -\int K(x - y) \cdot (\nabla g(y)) f(y) dy - \int (\nabla \cdot K(x - y)) g(y) f(y) dy.$$

\noindent Let $$G(x) = -\int K(x - y) \cdot (\nabla g(y)) f(y) dy$$ and $$H(x) = -\int (\nabla \cdot K(x - y)) g(y) f(y) dy.$$

\noindent Apply Holder's inequality and $|K(x - y)| \lesssim_{d} 1$ to $G(x)$ and the Hardy-Littlewood-Sobolev inequality, $|\nabla K(x - y)| \lesssim_{d} \frac{1}{|x - y|}$, and the Sobolev embedding theorem to $H(x)$. $\Box$\vspace{5mm}

\begin{corollary}\label{c4.5.6}
 $$(\ref{4.15.3}) \lesssim o(K).$$
\end{corollary}

\noindent \emph{Proof:} Let $$F(t,y) = \int_{\mathbf{R}^{d}} \frac{(x - y)_{j}}{|x - y|} \partial_{j} (e^{-ix \cdot \xi(t)} [\I(|u|^{4/d} u)(t,x)$$ $$- |\I u|^{4/d} (\I u)(t,x)]) (e^{ix \cdot \xi(t)} \overline{\I u}(t,x)) dx.$$

\noindent Then by lemma $\ref{l4.5.5}$,

\begin{equation}\label{4.18.1}
 F(t,y) = G(t,y) + H(t,y),
\end{equation}

\noindent with

\begin{equation}\label{4.18.2}
\aligned
\| G(t,y) \|_{L_{t}^{1} L_{x}^{\infty}([0, T] \times \mathbf{R}^{d})}  \lesssim \| \I(|u|^{4/d} u) - |\I u|^{4/d} (\I u) \|_{L_{t}^{2} L_{x}^{\frac{2d}{d + 2}}([0, T] \times \mathbf{R}^{d})} \\ \times \| \nabla e^{ix \cdot \xi(t)} \I u \|_{L_{t}^{2} L_{x}^{\frac{2d}{d - 2}}([0, T] \times \mathbf{R}^{d})} \lesssim o(K),
\endaligned
\end{equation}

\noindent and

\begin{equation}\label{4.18.3}
\aligned
\| H(t,y) \|_{L_{t}^{4/3} L_{x}^{\frac{6d}{5}}([0, T] \times \mathbf{R}^{d})}  \lesssim \| \I(|u|^{4/d} u) - |\I u|^{4/d} (\I u) \|_{L_{t}^{2} L_{x}^{\frac{2d}{d + 2}}([0, T] \times \mathbf{R}^{d})} \\ \| |\nabla|^{2/3} e^{ix \cdot \xi(t)} \I u \|_{L_{t}^{4} L_{x}^{\frac{2d}{d - 1}}([0, T] \times \mathbf{R}^{d})} \lesssim o(K^{2/3}).
\endaligned
\end{equation}

\noindent By Holder's inequality, conservation of mass,

$$\int_{0}^{T} \int_{\mathbf{R}^{d}} |\I u(t,y)|^{2} |G(t,y)| dy dt \leq \| G(t,y) \|_{L_{t}^{1} L_{x}^{\infty}([0, T] \times \mathbf{R}^{d})} \| \I u(t,y) \|_{L_{t}^{\infty} L_{x}^{2}([0, T] \times \mathbf{R}^{d})}^{2} \lesssim o(K).$$

\noindent By Sobolev embedding, lemma $\ref{l4.5.1}$,

$$\int_{0}^{T} \int_{\mathbf{R}^{d}} |\I u(t,y)|^{2} |H(t,y)| dy dt$$ $$\leq \| H(t,y) \|_{L_{t}^{4/3} L_{x}^{\frac{6d}{5}}([0, T] \times \mathbf{R}^{d})} \| \I u(t,y) \|_{L_{t}^{8} L_{x}^{\frac{12d}{6d - 5}}([0, T] \times \mathbf{R}^{d})}^{2} \lesssim o(K).$$

\noindent This implies $(\ref{4.15.3}) \lesssim o(K)$. $\Box$\vspace{5mm}

\noindent Finally consider $(\ref{4.15.4})$.

$$\I(|u|^{4/d} u) \overline{\I u} = |u|^{2 + 4/d} + (\I - 1)(|u|^{4/d} u)(\overline{\I u})$$ $$+ (1 - \I)(|u|^{4/d} u)\overline{(\I - 1)u} + \I(|u|^{4/d} u)\overline{(\I - 1)u}.$$

$$Im[|u|^{2 + 4/d}] \equiv 0.$$

\noindent Next, let

\begin{equation}\label{4.19.1}
 F_{1j}(t,x) = \int_{\mathbf{R}^{d}} \frac{(x - y)_{j}}{|x - y|} (1 - \I)(|u|^{4/d} u)(t,y)\overline{(\I - 1)u}(t,y) dy,
\end{equation}

\begin{equation}\label{4.19.2}
 F_{2j}(t,x) = \int_{\mathbf{R}^{d}} \frac{(x - y)_{j}}{|x - y|}  (\I - 1)(|u|^{4/d} u)(t,y)(\overline{\I u})(t,y) dy,
\end{equation}

\noindent and

\begin{equation}\label{4.19.3}
F_{3j}(t,x) = \int_{\mathbf{R}^{d}} \frac{(x - y)_{j}}{|x - y|} \I(|u|^{4/d} u)(t,y)\overline{(\I - 1)u}(t,y) dy.
\end{equation}

\noindent By Holder's inequality, lemma $\ref{l4.5.1}$, lemma $\ref{l4.5.2}$, and $|\xi - \xi(t)| \sim |\xi|$ for $|\xi| \geq CK$,

$$\| F_{1j} \|_{L_{t}^{1} L_{x}^{\infty}([0, T] \times \mathbf{R}^{d}}$$ $$\lesssim \| (1 - \I)u \|_{L_{t}^{2} L_{x}^{\frac{2d}{d - 2}}([0, T] \times \mathbf{R}^{d})} \| (1 - \I)(|u|^{4/d} u) \|_{L_{t}^{2} L_{x}^{\frac{2d}{d + 2}}([0, T] \times \mathbf{R}^{d})} \lesssim o(1).$$

\noindent Next, by lemma $\ref{l4.5.5}$, $e^{-ix \cdot \xi(t)}(\I - 1)(|u|^{4/d} u) = \nabla \cdot \frac{\nabla}{\Delta} e^{-ix \cdot \xi(t)} (\I - 1)(|u|^{4/d} u)$, we have $F_{2j} = G_{2j} + H_{2j}$, where

$$\| G_{2j} \|_{L_{t}^{1} L_{x}^{\infty}([0, T] \times \mathbf{R}^{d})}$$ $$\lesssim \| \nabla e^{-ix \cdot \xi(t)} \I u \|_{L_{t}^{2} L_{x}^{\frac{2d}{d - 2}}([0, T] \times \mathbf{R}^{d})} \| \frac{\nabla}{\Delta} e^{-ix \cdot \xi(t)} (\I - 1)(|u|^{4/d} u) \|_{L_{t}^{2} L_{x}^{\frac{2d}{d + 2}}([0, T] \times \mathbf{R}^{d})}$$

$$\lesssim o(K) (\frac{1}{K}) = o(1).$$

\noindent We use the fact that $|\xi - \xi(t)| \gtrsim K$ on the support of $(1 - \I)$.

$$\| H_{2j} \|_{L_{t}^{4/3} L_{x}^{\frac{6d}{5}}([0, T] \times \mathbf{R}^{d})}$$ $$\lesssim\| |\nabla|^{2/3} e^{-ix \cdot \xi(t)} \I u \|_{L_{t}^{2} L_{x}^{\frac{2d}{d - 2}}([0, T] \times \mathbf{R}^{d})} \| \frac{\nabla}{\Delta} e^{-ix \cdot \xi(t)} (\I - 1)(|u|^{4/d} u) \|_{L_{t}^{2} L_{x}^{\frac{2d}{d + 2}}([0, T] \times \mathbf{R}^{d})} \lesssim o(K^{-1/3}).$$

\noindent Finally, $F_{3j} = G_{3j} + H_{3j}$, with

$$\| G_{3j} \|_{L_{t}^{1} L_{x}^{\infty}([0, T] \times \mathbf{R}^{d})}$$ $$\lesssim \|  \nabla e^{-ix \cdot \xi(t)} \I u \|_{L_{t}^{2} L_{x}^{\frac{2d}{d - 2}}([0, T] \times \mathbf{R}^{d})} \| \frac{\nabla}{\Delta} e^{-ix \cdot \xi(t)} (\I - 1)(|u|^{4/d} u) \|_{L_{t}^{2} L_{x}^{\frac{2d}{d + 2}}([0, T] \times \mathbf{R}^{d})} \lesssim o(1),$$

$$\| H_{3j} \|_{L_{t}^{1} L_{x}^{3d}([0, T] \times \mathbf{R}^{d})}$$ $$\lesssim \| |\nabla|^{2/3} e^{-ix \cdot \xi(t)} \I u \|_{L_{t}^{2} L_{x}^{\frac{2d}{d - 2}}([0, T] \times \mathbf{R}^{d})} \| \frac{\nabla}{\Delta} e^{-ix \cdot \xi(t)} (\I - 1)(|u|^{4/d} u) \|_{L_{t}^{2} L_{x}^{\frac{2d}{d + 2}}([0, T] \times \mathbf{R}^{d})} \lesssim o(K^{-1/3}).$$

\noindent By Holder's inequality,

$$\int_{0}^{T} \int_{\mathbf{R}^{d}} |F_{1j}(t,x) + G_{2j}(t,x) + G_{3j}(t,x)| |(\nabla - i \xi(t)) \I u(t,x)| |\I u(t,x)| dx dt$$

$$\lesssim \| F_{1j}(t,x) + G_{2j}(t,x) + G_{3j}(t,x) \|_{L_{t}^{1} L_{x}^{\infty}([0, T] \times \mathbf{R}^{d})}$$ $$\times \| \nabla e^{-ix \cdot \xi(t)} \I u \|_{L_{t}^{\infty} L_{x}^{2}([0, T] \times \mathbf{R}^{d})} \| \I u(t,x) \|_{L_{t}^{\infty} L_{x}^{2}([0, T] \times \mathbf{R}^{d})} \lesssim o(K).$$

\noindent Next, by Holder's inequality and Sobolev embedding,

$$\int_{0}^{T} \int_{\mathbf{R}^{d}} |H_{3j}(t,x)| |\nabla e^{-ix \cdot \xi(t)} \I u(t,x)| |\I u(t,x)| dx dt$$ $$\lesssim \| H_{3j} \|_{L_{t}^{1} L_{x}^{3d}([0, T] \times \mathbf{R}^{d})} \| \nabla e^{-ix \cdot \xi(t)} \I u(t,x) \|_{L_{t}^{\infty} L_{x}^{2}([0, T] \times \mathbf{R}^{d})} \| e^{-ix \cdot \xi(t)} \I u \|_{L_{t}^{\infty} L_{x}^{\frac{6d}{3d - 2}}([0, T] \times \mathbf{R}^{d})}$$

$$\lesssim o(K^{-1/3}) K K^{1/3} = o(K).$$

\noindent Finally, by the Sobolev embedding theorem, lemma $\ref{l4.5.1}$, and interpolation,

$$\int_{0}^{T} \int_{\mathbf{R}^{d}} |H_{2j}(t,x)| |\nabla e^{-ix \cdot \xi(t)} \I u(t,x)| |\I u(t,x)| dx dt$$

$$\lesssim \| H_{2j} \|_{L_{t}^{4/3} L_{x}^{\frac{6d}{5}}([0, T] \times \mathbf{R}^{d})} \| \nabla e^{-ix \cdot \xi(t)} \I u \|_{L_{t}^{4} L_{x}^{\frac{2d}{d - 1}}([0, T] \times \mathbf{R}^{d})} \| \I u \|_{L_{t}^{\infty} L_{x}^{\frac{6d}{3d - 2}}([0, T] \times \mathbf{R}^{d})}$$

$$\lesssim o(K^{-1/3}) K K^{1/3} = o(K).$$

\noindent This completes the proof of theorem $\ref{t4.3}$. $\Box$\vspace{5mm}

\noindent Therefore, scenario $\int_{0}^{\infty} N(t)^{3} dt = \infty$ has been excluded.\vspace{5mm}

\noindent We have actually proved a more general estimate.

\begin{theorem}\label{t4.6}
 Suppose $a_{j}(t,x)$ is an odd function on $\mathbf{R}^{d}$ for all $t$ and there exists a constant $C$ such that

\begin{equation}\label{4.20.1}
 |a_{j}(t,x)| \leq C,
\end{equation}

\begin{equation}\label{4.20.2}
 |\partial_{k} a_{j}(t,x)| \leq \frac{C}{|x|}.
\end{equation}

\noindent Suppose also that $u(t,x)$ is a minimal mass blowup solution to $(\ref{0.1})$, $\mu = \pm 1$. Then

\begin{equation}\label{4.21.1}
\aligned
\int_{0}^{T} \int_{\mathbf{R}^{d} \times \mathbf{R}^{d}} a_{j}(t, x - y) |\I u(t,y)|^{2} Re \{ [\I(|u|^{4/d} \bar{u})(t,x) - |\I u|^{4/d} (\overline{\I u})(t,x)] \\ \times (\partial_{j} - i \xi_{j}(t)) \I u(t,x) \} dx dy dt \lesssim_{m_{0}, d} o(K) C,
\endaligned
\end{equation}

\begin{equation}\label{4.15.3}
\aligned
\int_{0}^{T} \int_{\mathbf{R}^{d} \times \mathbf{R}^{d}} a_{j}(t, x - y) |\I u(t,y)|^{2} Re \{ \overline{\I u}(t,x) (\partial_{j} - i \xi_{j}(t))\\ \times [|\I u|^{4/d} (\I u)(t,x) - \I(|u|^{4/d} u)(t,x)] \} dx dy dt \lesssim_{m_{0}, d} o(K) C,
\endaligned
\end{equation}

\begin{equation}\label{4.15.4}
\aligned
\int_{0}^{T} &\int_{\mathbf{R}^{d} \times \mathbf{R}^{d}} a_{j}(t, x - y) Re[[\overline{\I u}(t,x) (\partial_{j} - i \xi_{j}(t)) \I u(t,x)] \\ &[\I(|u|^{4/d} \bar{u})(t,y) \I u(t,y)  - \I(|u|^{4/d} u)(t,y) \overline{\I u}(t,y)]] dx dy dt \lesssim_{m_{0}, d} o(K) C.
\endaligned
\end{equation}

\end{theorem}

\noindent \emph{Proof:} By theorem $\ref{t2.3}$ a minimal mass blowup solution to $(\ref{0.1})$, $\mu = \pm 1$ satisfies the hypotheses of lemma $\ref{l4.5.1}$. $\Box$\vspace{5mm}

\noindent \textbf{Remark:} We conclude this section with a brief summary of what we have done. We have excluded the scenario when $\mu = +1$, $\int_{0}^{\infty} N(t)^{3} dt = \infty$ by proving that the errors arising from the interaction Morawetz estimates $(\ref{4.0.0})$, $(\ref{4.0.0.1})$ are bounded by $o(K)$. In the defocusing case these interaction Morawetz estimates are positive definite and $\gtrsim K$, which is a contradiction for $K$ sufficiently large. In the focusing case $(\ref{4.0.0})$ and $(\ref{4.0.0.1})$ are not positive definite. However, theorem $\ref{t4.6}$ states that if we did find an appropriate interaction Morawetz potential that satisfies $(\ref{4.20.1})$, $(\ref{4.20.2})$, then the error would be bounded by $o(K)$.


\section{$\int_{0}^{\infty} N(t)^{3} dt < \infty$} In this section we exclude the existence of a minimal mass blowup solution with $N(t) \leq 1$ and

\begin{equation}\label{3.0.0}
\int_{0}^{\infty} N(t)^{3} dt = K < \infty.
 \end{equation}

\noindent Excluding this scenario concludes the proof of theorem $\ref{t0.2}$. As in \cite{TVZ2} and \cite{KVZ} we will prove additional regularity. Conservation of energy precludes $N(t) \rightarrow 0$ as $t \rightarrow \infty$, giving a contradiction. To that end we prove:\vspace{5mm}

\begin{theorem}\label{t3.1}
Suppose $\int_{0}^{\infty} N(t)^{3} dt = K < \infty$, $\xi(0) = 0$, and $u$ is a minimal mass blowup solution to $(\ref{0.1})$, $\mu = \pm 1$. Then $u(t,x) \in H_{x}^{s}(\mathbf{R}^{d})$ for $0 \leq s < 1 + 4/d$ and $$\| u(t,x) \|_{L_{t}^{\infty} H_{x}^{s}((0, \infty) \times \mathbf{R}^{d})} \lesssim K^{s+}.$$
\end{theorem}

\noindent Recall from $(\ref{1.8})$ that we also have $$\sum_{J_{k}} N(J_{k}) \sim K.$$ This implies $|\xi(t_{1}) - \xi(t_{2})| \lesssim_{d} K$ for all $t_{1}, t_{2} \in (0, \infty)$. Therefore $|\xi(t)| \lesssim_{d} K$ on $(0, \infty)$.

\begin{theorem}\label{t3.2}
If theorem $\ref{t3.1}$ is true, a minimal mass blowup solution to $(\ref{0.1})$, $\mu = +1$, with $N(t) \leq 1$ and $$\int_{0}^{\infty} N(t)^{3} dt = K < \infty$$ does not exist.
\end{theorem}

\noindent \emph{Proof:} Recall the compactness modulus function $C(\eta)$ defined for all $0 < \eta < \infty$ from $(\ref{6.2})$ and $(\ref{6.3})$. There exists a function $\eta(t)$ such that for $1 < s < 1 + 4/d$,

\begin{equation}\label{3.0.1}
\lim_{t \rightarrow \pm \infty} C(\eta(t)) N(t) + \eta(t)^{\frac{s - 1}{2s}} = 0.
\end{equation}

\noindent So for any $\delta > 0$ there exists $T$ sufficiently large so that $$C(\eta(T)) N(T) + \eta(T)^{\frac{s - 1}{2s}} < \delta.$$ Make a Galilean transformation setting $\xi(T) = 0$.

\begin{equation}\label{3.0.2}
\aligned
\| u(T) \|_{\dot{H}^{1}(\mathbf{R}^{d})} \lesssim \| u_{|\xi| \leq C(\eta(T)) N(T)} \|_{\dot{H}^{1}(\mathbf{R}^{d})} + \| u_{|\xi| \geq C(\eta(T)) N(T)} \|_{\dot{H}^{1}(\mathbf{R}^{d})} \lesssim C(\eta(T)) N(T) + \eta(T)^{\frac{s - 1}{2s}}.
\endaligned
\end{equation}

\noindent The estimate on $u_{|\xi| \geq C(\eta(T)) N(T)}$ follows from interpolating $\| u_{|\xi| \geq C(\eta(T)) N(T)} \|_{L^{2}(\mathbf{R}^{d})} < \eta(T)^{1/2}$ with

\begin{equation}\label{3.0.3}
\| u(t) \|_{L_{t}^{\infty} \dot{H}_{x}^{s}((0, \infty) \times \mathbf{R}^{d})} \lesssim K^{s+}
\end{equation}

\noindent for $1 < s < 1 + 4/d$. Before we made the Galilean transformation that set $\xi(T) = 0$, we had $|\xi(t)| \lesssim_{d} K$ for all $t \in (0, \infty)$, so by the triangle inequality and $(\ref{3.0.3})$, after the Galilean transformation,

\begin{equation}\label{3.0.3.1}
\| u(T) \|_{\dot{H}_{x}^{s}(\mathbf{R}^{d})} \lesssim K^{s+}.
\end{equation}

\noindent This bound is uniform for $T \in [0, \infty)$. Also, by the Sobolev embedding theorem,

\begin{equation}\label{3.0.4}
\| u(T) \|_{L_{x}^{\frac{2(d + 2)}{d}}(\mathbf{R}^{d})}^{\frac{2(d + 2)}{d}} \lesssim_{d} \| u(T) \|_{\dot{H}_{x}^{\frac{d}{d + 2}}(\mathbf{R}^{d})}^{\frac{2(d + 2)}{d}} \lesssim_{d} \| u(T) \|_{\dot{H}^{1}(\mathbf{R}^{d})}^{2} \| u(T) \|_{L^{2}(\mathbf{R}^{d})}^{4/d} \lesssim_{m_{0}, d} \delta^{2}.
\end{equation}

\noindent Using conservation of energy and $(\ref{0.2.1})$, for all $t \in (0, \infty)$,

\begin{equation}\label{3.0.5}
E(u(T)) = E(u(t)) = \frac{1}{2} \int |\nabla u(t,x)|^{2} dx + \frac{d}{2d + 4} \int |u(t,x)|^{\frac{2d + 4}{d}} dx \lesssim_{m_{0},d} \delta^{2}.
\end{equation}

\noindent By $(\ref{6.2})$ and conservation of mass,

$$\frac{99 m_{0}^{2}}{100} < \int_{|x - x(0)| < \frac{1}{N(0)} C(\frac{m_{0}^{2}}{100})} |u(0,x)|^{2} dx,$$

\noindent which by Holder's inequality and conservation of energy,

$$\leq \frac{1}{N(0)^{\frac{2d}{d + 2}}} C(\frac{m_{0}^{2}}{100})^{\frac{2d}{d + 2}} \| u(0) \|_{L_{x}^{\frac{2(d + 2)}{d}}(\mathbf{R}^{d})}^{2} \leq \frac{1}{N(0)^{\frac{2d}{d + 2}}} C(\frac{m_{0}^{2}}{100})^{\frac{2d}{d + 2}} E(T)^{\frac{d}{d + 2}} \lesssim \frac{1}{N(0)^{\frac{2d}{d + 2}}} \cdot C(\frac{m_{0}^{2}}{100})^{\frac{2d}{d + 2}} \delta^{\frac{d}{(d + 2)}}.$$

\noindent For $\delta > 0$ very small this is a contradiction. Therefore theorem $\ref{t3.2}$ has been proved, assuming theorem $\ref{t3.1}$ is true. $\Box$\vspace{5mm}

\noindent \textbf{Remark:} We could also apply this argument to the case $\mu = -1$, $\| u_{0} \|_{L^{2}(\mathbf{R}^{d})}$ is less than the mass of the ground state. We will not bother to do that here.\vspace{5mm}

\noindent \emph{Proof of theorem $\ref{t3.1}$:} We will rely on two intermediate lemmas to prove theorem $\ref{t3.2}$. As usual we will partition $(0, \infty)$ into subintervals $J_{k}$ with $\| u \|_{L_{t,x}^{\frac{2(d + 2)}{d}}(J_{k} \times \mathbf{R}^{d})} = \epsilon$.\vspace{5mm}

\begin{lemma}\label{l3.2}
For any $1/2 \leq \lambda < 1 + 4/d$ and $\lambda \leq 1/2 + \sigma$, if

\begin{equation}\label{3.6}
\sup_{J_{k}} \| u_{> M} \|_{S_{\ast}^{0}(J_{k} \times \mathbf{R}^{d})} \lesssim_{m_{0}, d, \sigma} \frac{K^{\sigma}}{M^{\sigma}},
\end{equation}

\noindent then

\begin{equation}\label{3.7}
\| P_{|\xi| \geq N} (|u|^{4/d} u) \|_{S^{0}((0, \infty) \times \mathbf{R}^{d})} \lesssim_{m_{0}, d, \lambda} \frac{K^{\lambda}}{M^{\lambda}}.
\end{equation}
\end{lemma}

\noindent \emph{Proof:} We have already proved that for any compact interval $J$, when $N \leq K$,

\begin{equation}\label{3.1}
\| P_{|\xi - \xi(t)| > N} u \|_{L_{t}^{2} L_{x}^{\frac{2d}{d - 2}}(J \times \mathbf{R}^{d})} \lesssim \frac{K^{1/2}}{N^{1/2}}.
\end{equation}

\noindent Let $C$ be a large, fixed constant such that $|\xi(t_{1}) - \xi(t_{2})| \leq \frac{C}{1000} K$. When $N \leq CK$, take $J_{n} = [0, T_{n}]$. By theorem $\ref{t2.3}$, with implied constant independent of $T_{n}$,

\begin{equation}\label{3.1.1}
\| u_{|\xi - \xi(t)| \geq N} \|_{L_{t}^{2} L_{x}^{\frac{2d}{d - 2}}([0, T_{n}] \times \mathbf{R}^{d})} \lesssim_{m_{0}, d} \frac{K^{1/2}}{N^{1/2}}.
\end{equation}

\noindent Taking $T_{n} \rightarrow \infty$, we have

\begin{equation}\label{3.1.1.1}
\| u_{|\xi - \xi(t)| \geq N} \|_{L_{t}^{2} L_{x}^{\frac{2d}{d - 2}}((0, \infty) \times \mathbf{R}^{d})} \lesssim_{m_{0},d} \frac{K^{1/2}}{N^{1/2}}.
\end{equation}

\noindent In fact, for any $\lambda \geq 1/2$, when $N \leq CK$,

\begin{equation}\label{3.1.2}
\| u_{|\xi - \xi(t)| \geq N} \|_{L_{t}^{2} L_{x}^{\frac{2d}{d - 2}}((0, \infty) \times \mathbf{R}^{d})} \lesssim_{m_{0}, d, \lambda} \frac{K^{\lambda}}{N^{\lambda}}.
\end{equation}

\noindent Interpolating this with conservation of mass,

\begin{equation}\label{3.1.3}
\| u_{|\xi - \xi(t)| \geq N} \|_{S^{0}((0, \infty) \times \mathbf{R}^{d})} \lesssim_{m_{0}, d, \lambda} \frac{K^{\lambda}}{N^{\lambda}}.
\end{equation}

\noindent Now we can use

\begin{lemma}\label{l3.2.1}
 Let $u$ be a solution to $(\ref{0.1})$ which is almost periodic modulo scaling on its maximal lifespan $I$, $u$ blows up forward in time. Then for all $t \in I$,

\begin{equation}\label{3.1.4}
 u(t) = \lim_{T \nearrow \sup I} i \int_{t}^{T} e^{i(t - \tau) \Delta} F(u(\tau)) d\tau,
\end{equation}

\noindent as a weak limit in $L_{x}^{2}$. 

\end{lemma}

\noindent \emph{Proof:} See section 6 of \cite{TVZ1}. $\Box$\vspace{5mm}

\noindent For $N \geq CK$,

\begin{equation}\label{3.1.5}
\| P_{|\xi - \xi(t)| > N} u \|_{S^{0}((0, \infty) \times \mathbf{R}^{d})} \lesssim_{d} \| P_{|\xi| \geq \frac{N}{2}} (|u|^{4/d} u) \|_{L_{t}^{2} L_{x}^{\frac{2d}{d + 2}}([0, \infty) \times \mathbf{R}^{d})}.
\end{equation}

\begin{equation}\label{3.4}
\aligned
&\lesssim \| P_{|\xi| \geq \frac{N}{2}} (|u_{|\xi - \xi(t)| \leq \eta N}|^{4/d} u_{|\xi - \xi(t)| \leq \eta N}) \|_{L_{t}^{2} L_{x}^{\frac{2d}{d + 2}}((0,\infty) \times \mathbf{R}^{d})} \\
&+ \| (1 - \chi(t)) u \|_{L_{t}^{\infty} L_{x}^{2}((0, \infty) \times \mathbf{R}^{d})}^{4/d} \| u_{|\xi - \xi(t)| \geq \eta N} \|_{L_{t}^{2} L_{x}^{\frac{2d}{d - 2}}((0,\infty) \times \mathbf{R}^{d})} \\ &+ \| u_{|\xi - \xi(t)| \geq C_{0} N(t)} \|_{L_{t}^{\infty} L_{x}^{2}((0, \infty) \times \mathbf{R}^{d})}^{4/d} \| u_{|\xi - \xi(t)| \geq \eta N} \|_{L_{t}^{2} L_{x}^{\frac{2d}{d - 2}}((0,\infty) \times \mathbf{R}^{d})} \\
&+ \| (u_{|\xi - \xi(t)| \geq \eta N}) |\chi(t) u_{|\xi - \xi(t)| \leq C_{0} N(t)}|^{4/d} \|_{L_{t}^{2} L_{x}^{\frac{2d}{d + 2}}((0,\infty) \times \mathbf{R}^{d})}.
\endaligned
\end{equation}

\noindent Therefore, for any $0 < s < 1 + 4/d$,

\begin{equation}\label{3.5}
\aligned
(\ref{3.4}) \lesssim_{m_{0},d,s} \sum_{M \leq \eta N} (\frac{M}{N})^{s} \| u_{|\xi - \xi(t)| \geq M} \|_{L_{t}^{2} L_{x}^{\frac{2d}{d - 2}}((0, \infty) \times \mathbf{R}^{d})}
+ \delta(C_{0}) \| P_{|\xi - \xi(t)| \geq \eta N} \|_{L_{t}^{2} L_{x}^{\frac{2d}{d - 2}}((0, \infty) \times \mathbf{R}^{d})}
\endaligned
\end{equation}

\begin{equation}\label{3.5.1}
 \left\{
    \begin{array}{ll}
     + C_{0}^{3/2} (\frac{K}{\eta N})^{1/2} (\sup_{J_{k}} \| u_{|\xi - \xi(t)| \geq \frac{\eta N}{2}} \|_{S_{\ast}^{0}(J_{k} \times \mathbf{R}^{d})}), & \hbox{if $d = 3$;} \\
      + C_{0}^{4 - 6/d} (\frac{K}{\eta N})^{2/d} (\sup_{J_{k}} \| u_{|\xi - \xi(t)| \geq \frac{\eta N}{2}} \|_{S_{\ast}^{0}(J_{k} \times \mathbf{R}^{d})})^{4/d}  \| u_{|\xi - \xi(t)| \geq \eta N} \|_{L_{t}^{2} L_{x}^{\frac{2d}{d - 2}}((0, \infty) \times \mathbf{R}^{d})}^{1 - 4/d} , & \hbox{if $d \geq 4$.}
    \end{array}
  \right.
\end{equation}

\noindent By induction,

 $$\| P_{|\xi| \geq \frac{N}{2}} (|u|^{4/d} u) \|_{S^{0}((0, \infty) \times \mathbf{R}^{d})} \leq \sum_{M \leq \eta N} C_{2}(m_{0}, d,s) C_{3}(m_{0}, d, \lambda) (\frac{K}{M})^{\lambda} \eta^{s - \lambda}$$ $$+ \delta(C_{0}) C_{2}(m_{0},d,s) C_{3}(m_{0}, d, \lambda) (\frac{K}{\eta N})^{\lambda}$$ $$+ \left\{
                                                                   \begin{array}{ll}
                                                                     C_{2}(m_{0},d,s) (\frac{K^{1/2}}{N^{1/2}}) (\frac{K^{\sigma}}{N^{\sigma}}) \frac{C_{0}^{3/2}}{\eta^{1/2 + \sigma}}, & \hbox{if $d = 3$;} \\
                                                                     C_{2}(m_{0},d,s) C_{3}(m_{0}, d, \lambda)^{1 - 4/d} (\frac{K}{\eta N})^{2/d} (\frac{K}{\eta N})^{4 \sigma/d} (\frac{K}{\eta N})^{(1 - 4/d) \lambda} C_{0}^{4 - 6/d}, & \hbox{if $d \geq 4$.}
                                                                   \end{array}
                                                                 \right.
.$$

\noindent If $\lambda < 1 + 4/d$ we can find $s$ such that $\lambda < s < 1 + 4/d$. Take $s = \frac{\lambda + 1 + 4/d}{2}$. Choose $\eta$ sufficiently small so that $\eta^{s - \lambda} C_{2}$ is very small. Then take $C_{0}(d, s, \eta, \lambda)$ sufficiently large so that $\frac{\delta(C_{0})}{\eta^{\lambda}} C_{2}$ is very small. Finally, if $d = 3$ choose $C_{3}$ sufficiently large so that $$\frac{C_{2} C_{0}^{3/2}}{\eta^{\lambda}} << C_{3},$$ and if $d \geq 4$ choose $C_{3}$ sufficiently large so that $$\frac{C_{2} C_{0}^{4 - 6/d}}{\eta^{\lambda}} << C_{3}^{4/d}.$$

\noindent This closes the induction and completes the proof. $\Box$\vspace{5mm}

\noindent \textbf{Remark:} We assume $K << \eta N$, otherwise we just use the results of $\S 3$.\vspace{5mm}

\noindent Now suppose $I$ is some interval $[a, b]$.

\begin{lemma}\label{l3.4}
If $u$ is a solution to $(\ref{0.1})$, 

\begin{equation}\label{3.11.2.1}
\| P_{> N} u(a) \|_{L_{x}^{2}(\mathbf{R}^{d})} \lesssim_{m_{0},d, \lambda} \frac{K^{\lambda}}{N^{\lambda}},
\end{equation}

\noindent with $\lambda < 1 + 4/d$, and 

\begin{equation}\label{3.11.3}
\| u \|_{L_{t,x}^{\frac{2(d + 2)}{d}}(I \times \mathbf{R}^{d})} \leq \delta
\end{equation}

\noindent for some $\delta(m_{0}, d, \lambda) > 0$ sufficiently small, then

\begin{equation}\label{3.12}
\| P_{> N} u \|_{S_{\ast}^{0}(I \times \mathbf{R}^{d})} \lesssim_{m_{0}, d, \lambda} \frac{K^{\lambda}}{N^{\lambda}}.
\end{equation}
\end{lemma}

\noindent \emph{Proof:} By Duhamel's formula and $(\ref{5.6})$,

\begin{equation}\label{3.13}
\| P_{> N} u \|_{S_{\ast}^{0}(I \times \mathbf{R}^{d})} \equiv \| P_{> N} u(a) \|_{L_{x}^{2}(\mathbf{R}^{d})} + \| P_{> N} (|u|^{4/d} u) \|_{L_{t,x}^{\frac{2(d + 2)}{d + 4}}(I \times \mathbf{R}^{d})}.
\end{equation}

\noindent Since $$\| u \|_{S_{\ast}^{0}(J_{k} \times \mathbf{R}^{d})} \lesssim_{m_{0},d} 1 + \delta^{1 + 4/d} \lesssim_{m_{0},d} 1,$$ our lemma is true for $N \leq CK$. By Bernstein's inequality and corollary $\ref{c5.4.1}$,

\begin{equation}\label{3.14}
\aligned
\| P_{> N} (|u|^{4/d} u) \|_{L_{t,x}^{\frac{2(d + 2)}{d + 4}}(I \times \mathbf{R}^{d})} &\lesssim_{d} \| P_{> N}(|u_{\leq N}|^{4/d} u_{\leq N}) \|_{L_{t,x}^{\frac{2(d + 2)}{d + 4}}} + \| |u_{> N} | |u|^{4/d} \|_{L_{t,x}^{\frac{2(d + 2)}{d + 4}}}
\\ &\lesssim_{m_{0},d,s}  \sum_{M \leq N} (\frac{M}{N})^{s} \| P_{> M} u \|_{L_{t,x}^{\frac{2(d + 2)}{d}}(I \times \mathbf{R}^{d})} \| u \|_{L_{t,x}^{\frac{2(d + 2)}{d}}(I \times \mathbf{R}^{d})}^{4/d}
\\ &\lesssim_{m_{0},d,s} \sum_{M \leq N} (\frac{M}{N})^{s} \| P_{> M} u \|_{L_{t,x}^{\frac{2(d + 2)}{d}}(I \times \mathbf{R}^{d})} \delta^{4/d}.
\endaligned
\end{equation}

\noindent Then apply the method of continuity. Recursively define a sequence of functions,

\begin{equation}\label{3.14.0.0}
\aligned
u_{0} &= e^{it \Delta} u(0), \\
u_{n + 1} &= e^{it \Delta} u_{0} - i \int_{0}^{t} e^{i(t - \tau) \Delta} |u_{n}(\tau)|^{4/d} u_{n}(\tau) d\tau.
\endaligned
\end{equation}

\noindent By $(\ref{3.11.2.1})$ and Strichartz estimates,

$$\| P_{> N} e^{it \Delta} u(0) \|_{L_{t,x}^{\frac{2(d + 2)}{d}}(I \times \mathbf{R}^{d})} \lesssim_{d} \frac{K^{\lambda}}{N^{\lambda}}.$$ Let $s = \frac{\lambda + 1 + 4/d}{2}$,

$$\sup_{N}  (\frac{N}{K})^{\lambda} \| P_{> N} u_{n + 1} \|_{L_{t,x}^{\frac{2(d + 2)}{d}}(I \times \mathbf{R}^{d})} \lesssim_{m_{0},d,\lambda} 1+ \delta^{4/d} (\sup_{N} (\frac{N}{K})^{\lambda} \| P_{> N} u_{n} \|_{L_{t,x}^{\frac{2(d + 2)}{d}}(I \times \mathbf{R}^{d})} ).$$

\noindent By continuity, for $\delta(m_{0}, d, \lambda) > 0$ sufficiently small $$\| u_{> N} \|_{L_{t,x}^{\frac{2(d + 2)}{d}}(I \times \mathbf{R}^{d})} \lesssim_{m_{0},d,\lambda} \frac{K^{\lambda}}{N^{\lambda}}.$$ By the same argument we also have

\begin{equation}\label{3.14.0.0.1}
\| P_{> N} (|u|^{4/d} u) \|_{L_{t,x}^{\frac{2(d + 2)}{d + 4}}(I \times \mathbf{R}^{d})} \lesssim_{m_{0},d,\lambda} \frac{K^{\lambda}}{N^{\lambda}}.
\end{equation}

\noindent Therefore,

\begin{equation}\label{3.14.0.1}
\| u_{> N} \|_{S_{\ast}^{0}(I \times \mathbf{R}^{d})} \lesssim_{m_{0},d,\lambda} \frac{K^{\lambda}}{N^{\lambda}}.
\end{equation}

\noindent Again, since $|\xi| \sim |\xi - \xi(t)|$ when $N \geq CK$, this proves

\begin{equation}\label{3.14.0.2}
\| u_{|\xi - \xi(a)| > N} \|_{S_{\ast}^{0}(I \times \mathbf{R}^{d})} \lesssim \frac{K^{\lambda}}{N^{\lambda}}.
\end{equation}

$\Box$\vspace{5mm}

\begin{corollary}\label{c3.5}
If $u$ is a solution to $(\ref{0.1})$, $\lambda < 1 + 4/d$,

\begin{equation}\label{3.14.0.3}
\| P_{> N} u \|_{L_{t}^{\infty} L_{x}^{2}((0, \infty) \times \mathbf{R}^{d})} \lesssim \frac{K^{\lambda}}{N^{\lambda}},
 \end{equation}

 \noindent and

 \begin{equation}\label{3.14.0.4}
 \| u \|_{L_{t,x}^{\frac{2(d + 2)}{d}}(J_{k} \times \mathbf{R}^{d})} = \epsilon,
 \end{equation}

 \noindent then
 \begin{equation}\label{3.15.0.5}
 \| u_{> N} \|_{S_{\ast}^{0}(J_{k} \times \mathbf{R}^{d})} \lesssim \frac{K^{\lambda}}{N^{\lambda}}.
 \end{equation}
 \end{corollary}

\noindent \emph{Proof:} Partition each subinterval $J_{k}$ with $\| u \|_{L_{t,x}^{\frac{2(d + 2)}{d}}(J_{k} \times \mathbf{R}^{d})} = \epsilon$ into a finite number of subintervals $I_{i}$ with $\| u \|_{L_{t,x}^{\frac{2(d + 2)}{d}}(I_{i} \times \mathbf{R}^{d})} = \delta$. Combining $(\ref{3.14.0.3})$ and lemma $\ref{l3.4}$,

\begin{equation}\label{3.14.2}
\| P_{> N} u \|_{S_{\ast}^{0}(J_{k} \times \mathbf{R}^{d})} \lesssim \frac{K^{\lambda}}{N^{\lambda}}.
\end{equation}

\noindent Now we are ready to prove theorem $\ref{t3.1}$.\vspace{5mm}

\noindent \emph{Proof of theorem $\ref{t3.1}$:} This is proved by induction. Take $N \geq CK$. Lemma $\ref{l3.2}$ implies that since $\| u \|_{S_{\ast}^{0}(J_{k} \times \mathbf{R}^{d})} \lesssim 1$,

$$\| u_{|\xi - \xi(t)| > N} \|_{S^{0}((0, \infty) \times \mathbf{R}^{d})} \lesssim_{m_{0},d} \frac{K^{1/2}}{N^{1/2}}.$$

\noindent By corollary $\ref{c3.5}$ this implies

$$\| P_{|\xi - \xi(a_{k})| > N} u \|_{S_{\ast}^{0}(J_{k} \times \mathbf{R}^{d})} \lesssim \frac{K^{1/2}}{N^{1/2}}.$$

\noindent Applying lemma $\ref{l3.2}$ again we have $$\| u_{|\xi - \xi(t)| > N} \|_{S^{0}((0, \infty) \times \mathbf{R}^{d})} \lesssim \frac{K}{N}.$$ Iterating at most four more times, theorem $\ref{t3.1}$ is proved. $\Box$\vspace{5mm}

\noindent We have excluded the second minimal mass blowup scenario. This concludes the proof of theorem $\ref{t0.2}$. $\Box$

\newpage

\nocite*
\bibliographystyle{plain}
\bibliography{dgeq3}

\end{document}